%% file: FieldClassNoAvg_2015_V14.tex
\documentclass[oneside]{amsart}
\usepackage{amssymb,latexsym,amsmath,amsthm,enumerate,mathrsfs}
\usepackage{fancyhdr}
\usepackage{color}
\usepackage{stmaryrd}
\usepackage{fullpage}

\input{format}
\input{melanie.sty}

\newcommand{\pfrak}{\mathfrak{p}}
\newcommand{\Cl}{\mathrm{Cl}}

\newcommand{\Ocal}{\mathcal{O}}
\newcommand{\Ascr}{\mathscr{A}}
\newcommand{\Bscr}{\mathscr{B}}
\newcommand{\Pscr}{\mathscr{P}}
\newcommand{\Escr}{\mathscr{E}}

\newcommand{\Leg}[2]{\left( \frac{#1}{#2}\right)}

\definecolor{pink}{rgb}{1,.2,.6}
\definecolor{orange}{rgb}{0.7,0.3,0}
\definecolor{blue}{rgb}{.2,.6,.75}
\definecolor{green}{rgb}{.4,.7,.4}
\definecolor{purple}{RGB}{127,0,255}

\begin{document}

\numberwithin{equation}{section}

\title[On $\ell$-torsion in class groups]{On $\ell$-torsion in class groups\\ of number fields}

\author[Ellenberg]{Jordan Ellenberg}
\address{Department of Mathematics, 480 Lincoln Dr., Madison, WI 53706 USA}
\email{ellenber@math.wisc.edu}

\author[Pierce]{Lillian. B. Pierce}
\address{Department of Mathematics, Duke University, 120 Science Drive, Durham NC 27708}
\email{pierce@math.duke.edu}

\author[Wood]{Melanie Matchett Wood}
\address{Department of Mathematics, 480 Lincoln Dr., Madison, WI 53706 USA\\
and
American Institute of Mathematics\\600 East Brokaw Road\\
San Jose, CA 95112 USA}
\email{mmwood@math.wisc.edu}



\begin{abstract}
For each integer $\ell \geq 1$, we prove an unconditional upper bound on the size of the $\ell$-torsion subgroup of the class group, which holds for all but a zero-density set of field extensions of $\Q$ of degree $d$, for any fixed $d \in \{2,3,4,5\}$ (with the additional restriction in the case $d=4$ that the field be non-$D_4$).   For sufficiently large $\ell$ (specified explicitly), these results are as strong as a previously known bound that is conditional on GRH.  As part of our argument, we develop a probabilistic ``Chebyshev sieve,'' and give uniform, power-saving error terms for the asymptotics of quartic (non-$D_4$) and quintic fields with chosen splitting types at a finite set of primes.
 \end{abstract}

\maketitle

\section{Introduction}

The distribution of class groups is a great mystery.  The Cohen-Lenstra heuristics \cite{CohLen84} (for quadratic fields) and the Cohen-Lenstra-Martinet heuristics \cite{Cohen1990} (for more general number fields)
make predictions for the distribution of class groups, including for the average size of the $\ell$-torsion subgroups  for certain ``good'' primes $\ell$.
 However, the questions of proving anything towards these predictions are almost entirely open, and mostly apparently inaccessible.

The main goal of the present work is to prove, for each integer $\ell \geq 1$, an unconditional upper bound for the size of the $\ell$-torsion subgroup of the class group, which holds for all but a zero-density set of field extensions of $\Q$ of degree $d$, for any fixed $d \in \{2,3,4,5\}$ (with the additional restriction in the case $d=4$ that the field be non-$D_4$). 
Alternatively, these results may be viewed as the first unconditional upper bounds for the average size of $\ell$-torsion in class groups as the field varies over extensions of $\Q$ of fixed degree $d\in \{2,3,4,5\}$ (and non-$D_4$ in the case $d=4$).

Let $K$ be a degree $d$ field extension of $\Q$ with absolute discriminant $D_K = | \disc K/\Q|$. 
We will denote the class group by   $\Cl_K$ and the $\ell$-torsion subgroup by $\Cl_K[\ell]$.
We note the trivial pointwise upper bound (see for example \cite[Thm 4.4]{Nar80})
\beq\label{Cl_trivial}
|\Cl_K[\ell]| \leq |\Cl_K| \ll_{d,\ep} D_K^{1/2 + \ep},
\eeq
for every $\ep>0$.
 (Throughout, $\ep>0$ is allowed to be arbitrarily small (possibly taking a different value in different occurrences), and $A \ll B$ indicates that $|A| \leq cB$ for an implied constant $c$, which we allow in any instance to depend on $\ell, d, \ep$.)

It is conjectured that 
\beq\label{Cl_conj}
|\Cl_K[\ell]| \ll D_K^{\ep}
\eeq
for every $\ep>0$, but improving on the trivial bound (\ref{Cl_trivial}) has proved difficult.
(Impetus for this conjecture may be found in Duke \cite{Duk98}, Zhang \cite[page 10]{Zha05}, and Brumer and Silverman \cite[``Question $\mathrm{CL}(\ell,d)$'']{BruSil96}.) 
For $K$ quadratic, Gauss's genus theory \cite{Gau} implies (\ref{Cl_conj}) in the case $\ell=2$.
Recently, \cite{BSTTTZ17} obtained nontrivial upper bounds for 2-torsion in fields of degree $d$ for all $d \geq 3$, proving $|\Cl_K[2]|\ll D_K^{0.2784... + \ep}$ for $d=3,4$ and $|\Cl_K[2]|\ll D_K^{1/2 - 1/2d+ \ep}$ for $d \geq 5$.
For $\ell=3$, after initial incremental improvement  in \cite{HelVen06}, \cite{Pie05}, \cite{Pie06} over the trivial bound (\ref{Cl_trivial}) for quadratic fields, Ellenberg and Venkatesh proved \cite[Prop. 3.4, Cor. 3.7]{EllVen07} 
 that 
 \beq\label{cubic_3_del}
 | \Cl_K[3]| \ll D_K^{1/3+\ep}
 \eeq
holds for both quadratic and cubic fields, and moreover 
 there is a positive constant $\del>0$ such that 
 \beq\label{cubic_4_del}
 |\Cl_K[3]| \ll D_K^{1/2 - \del + \ep}
 \eeq
holds for quartic fields. (In particular, one may take $\del = 1/168$ in (\ref{cubic_4_del}) for quartic fields with Galois closure having Galois group $A_4$ or $S_4$.) At this time, these are the best  bounds in the literature that are unconditional and hold for all such fields.

In the realm of average results, there 
is little known, with the exceptions being spectacular successes. 
For 3-torsion in quadratic fields, Davenport and Heilbronn \cite{DavHei71} proved
\beq\label{DavHei71_3} \sum_{\bstack{\deg(K)=2}{0<D_K\leq X}} |\Cl_{K}[3]| \sim \left( \frac{2}{3\zeta(2)} + \frac{1}{\zeta(2)}\right) X,\eeq
in which the first contribution is from fields with $\disc K/\Q >0$ and the second is from fields with $\disc K/\Q <0$; this has recently been improved to reflect second order terms by \cite{BST}, \cite{TanTho13} and \cite{Hou13}.
For 2-torsion in cubic fields, Bhargava \cite{Bha05} proved the asymptotic
\beq\label{Bhargava_2} \sum_{\bstack{\deg(K)=3}{0<D_K\leq X}}| \Cl_{K}[2] |\sim \left( \frac{5}{48\zeta(3)} + \frac{3}{8\zeta(3)} \right) X,
\eeq
in which each isomorphism class of fields is counted once, and the first contribution is from fields with $\disc K/\Q >0$ and the second is from fields with $\disc K/\Q<0$. 
For $4$-torsion in quadratic fields, Fouvry and Kl\"{u}ners \cite{FK} have determined the asymptotics, for each non-negative integer $k$,
\beq\label{FK} \sum_{\bstack{\deg(K)=2}{0<D_K\leq X}} |\Cl_{K}^2[2]|^k \sim (c_k + p^{-k}(c_{k+1}-c_k) ) X,\eeq
where $c_k$ is the number of vector subspaces of $\F_2^k$.
See also the recent work of Klys \cite{Klys} giving analogous results on $3$-torsion in cyclic cubic fields, and the recent work of Milovic on 16-rank in quadratic fields, e.g. \cite{Mil15}.

Turning to conditional results,
 Klys's results \cite{Klys} extend to $p$-torsion in cyclic degree $p$ fields under GRH and Smith \cite{Smith} has results on $8$-torsion averages in quadratic fields under GRH as well.
 In the case of quadratic fields, Wong \cite{Won99} proved that, conditional on the Birch--Swinnerton-Dyer conjecture and the Riemann Hypothesis, 
$|\Cl_K[3]| \ll D_K^{1/4+\ep}.$
Before the proof of (\ref{cubic_3_del}), Soundararajan noted (as communicated in \cite{HelVen06}) that one could prove 
 $|\Cl_K[3]| \ll D_K^{1/3+\ep}$ for $K$ quadratic
if one assumed the truth of the Riemann Hypothesis  for only the $L$-function $L(s,\chi)$ of the quadratic character $\chi$ associated to the quadratic field $K$. The key idea of the latter bound was the use of  many small primes that split in $K$; the role of the Riemann Hypothesis was to guarantee the existence of sufficiently many such primes. This approach has been generalized by Ellenberg and Venkatesh  \cite{EllVen07} to number fields of any degree; we recall the key result in the special case of field extensions of $\Q$:

\begin{letterthm}[Lemma 2.3 of \cite{EllVen07}]\label{lemma_M_primes}
Let $K$ be a field extension of $\Q$ of degree $d$, and let $\ell$ be a positive integer. Let $\del< \frac{1}{2\ell(d-1)}$. Suppose that $\{\pfrak_1,\ldots, \pfrak_M\}$ are $M$ prime ideals in $\Ocal_K$ with $\mathrm{Norm}(\pfrak_j) \leq D_K^{\del}$ that are unramified in $K/\Q$ and are not extensions of ideals from any proper subfield $K_0 \subsetneq K$.Then 
\beq\label{EllVen_ineq}
 |\Cl_K[\ell]| \ll_{d,\ell,\ep} D_K^{1/2 + \ep}M^{-1}.
 \eeq
\end{letterthm}
(Here we recall the convention in \cite{EllVen07} that an ideal $\pfrak$ in $\Ocal_K$ is said to be an extension of a prime ideal from a subfield $K_0 \subsetneq K$ if there is a prime ideal $\pfrak_0$ in $\Ocal_{K_0}$ such that $\pfrak= \pfrak_{0}\Ocal_{K} $.)

Upon assuming GRH, an application of the effective Chebotarev theorem of Lagarias and Odlyzko \cite{LagOdl75} guarantees, for any fixed $\eta >0$, the existence of $\gg D_K^{\eta-\ep}$ rational primes of size $\leq D_K^{\eta}$ that split completely in $K$. Upon choosing $\eta  = \frac{1}{2\ell(d-1)} -\ep_0$ for arbitrarily small $\ep_0>0$, one obtains the following bound, currently  the state of the art for conditional  pointwise upper bounds for $|\Cl_K[\ell]|$:
\begin{letterthm}[Prop. 3.1 of  \cite{EllVen07}]\label{thm_EllVen07}
Let $K$ be a field extension of $\Q$ of degree $d$ and $\ell$ a positive integer. Assuming GRH,
\beq\label{EV_GRH}
|\Cl_K[\ell]| \ll_{d,\ell,\ep} D_K^{\frac{1}{2} - \frac{1}{2\ell(d-1)} + \ep},
\eeq
for any $\ep>0$.
\end{letterthm}

One may attempt to remove the conditionality by proving results that hold on average, or for all but a small exceptional family.
In this vein, in the case of imaginary quadratic fields, Soundararajan \cite{Sou00} noted that for all but at most one imaginary quadratic field $K$ with $D_K \in [X,2X]$, one has the bound $|\Cl_{K}[\ell]| \ll X^{1/2-1/2\ell + \ep}$, for any prime $\ell$.
Also in the imaginary quadratic case, 
a recent result of Heath-Brown and Pierce \cite{PieHB14a} provides an upper bound for averages (and in addition higher moments) of $|\Cl_K[\ell]|,$  for example proving for any prime $\ell \geq 5$ that
\beq\label{HB_Pie_result}
\sideset{}{'}\sum_{\bstack{\deg(K)=2}{0<D_K\leq X}} |\Cl_K[\ell]| \ll X^{\frac{3}{2} - \frac{3}{2\ell+2} + \ep},
\eeq
with the sum restricted to imaginary quadratic fields.

In this paper, we prove \emph{unconditional} results for $|\Cl_K[\ell]|$ that  are as strong as (\ref{EV_GRH}) for all sufficiently large positive integers $\ell$, and hold for all but a zero-density family of   quadratic, cubic, non-$D_4$-quartic, or quintic field extensions of $\Q$.

For this we work with families of fields. Let $N_d(X)$ denote the number of degree $d$ extensions of $\Q$ with $0<D_K\leq X$, in which each isomorphism class is counted once; it is conjectured that for an appropriate constant $c_d$,
\beq\label{num_fields}
 N_d(X) \sim c_{d}X.
 \eeq
Importantly for our work, this is known to be true for $d=2$ (classical), $d=3$ by Davenport and Heilbronn \cite{DavHei71},
$d=4$ by Cohen, Diaz y Diaz, and Olivier \cite{CDO02} and Bhargava \cite{Bha05},
 and  $d=5$ by Bhargava \cite{Bha10a}. 
Throughout our work, in the case of $d=4$, we restrict our attention to non-$D_4$-quartic fields (that is, quartic extensions whose Galois closure does not have Galois group $D_4$); see Remark~\ref{rem_d4}. Thus we let $\tilde{N}_4(X)$ denote the further restricted count of non-$D_4$-quartic extensions of $\Q$; then (\ref{num_fields}) is also known to hold for $\tilde{N}_4(X)$, with a different constant \cite{Bha05}. 

As a consequence of the field counts (\ref{num_fields}) combined with the trivial bound (\ref{Cl_trivial}), a trivial average bound for $|\Cl_K[\ell]|$ is
\beq\label{trivial_avg}
\sum_{\bstack{\deg(K)=d}{0 < D_K \leq X}}  |\Cl_K[\ell]|  \ll_{d,\ep} X^{3/2 + \ep}.
 \eeq
Our approach to improve upon (\ref{trivial_avg}) is to show that ``most'' degree $d$ fields $K$ contain sufficiently many small primes that split completely in $K$ for Theorem \ref{lemma_M_primes} to give a good upper bound for $|\Cl_K[\ell]|$. Roughly speaking, we will show that there is some small $\del_0>0$ such that  for all but at most $O(X^{1- \del_0})$ of the degree $d$ fields $K$ with $0<D_K \leq X $, at least  a fixed positive proportion of the primes $p \leq X^{\del_0}$ split completely in $K$.  (Under GRH, the small set of exceptional fields is in fact empty.)

Our main results are  as follows:
\begin{thm}\label{thm_main}
Let $d \in \{2,3,4,5\}$ and let $\ell$ be any positive integer with $\ell \geq \ell(d)$ where 
\[ \ell(2)=  \ell(3) =1, \qquad \ell(4) =8 , \qquad \ell(5)= 25.\]  Then for all but $O_{d,\ell,\ep}(X^{1 - \frac{1}{2\ell (d-1)}+\ep})$ degree $d$ fields $K/\Q$ with $D_K \leq X$ (and non-$D_4$ in the case $d=4$),
\[ |\Cl_K[\ell]| \ll_{d,\ell,\ep} D_K^{\frac{1}{2} - \frac{1}{2\ell(d-1)} + \ep},\]
for all $\ep>0$.
For $d=4,5,$ in the remaining cases of positive integers $\ell < \ell(d)$, for all but $O_{d,\ep}(X^{1 - \del_0(d)+\ep})$ degree $d$ fields $K/\Q$ with $0<D_K \leq X$ (and non-$D_4$ in the case $d=4$),
\[ |\Cl_K[\ell]| \ll_{d,\ell,\ep} D_K^{\frac{1}{2} - \del_0(d) + \ep},\]
for all $\ep>0$,
where we may take
\[
 \del_0(d) = \begin{cases}
	1/48 & \text{if $d=4$} \\
	1/200 & \text{if $d=5$}.
		\end{cases}
		\]
\end{thm}
\begin{remark}
Theorem \ref{thm_B} states a version of this result in terms of bounding the number of exceptional fields that fail to have many small split primes. One notes from Theorem \ref{thm_B} that for sufficiently large $\ell$, the limiting reagent is not the availability of small completely split primes, but the constraint $\del< (2\ell(d-1))^{-1}$ in Theorem \ref{lemma_M_primes}.
\end{remark}

As immediate corollaries, we note:
\begin{cor}\label{cor_main_d}
Let $d \in \{2,3,4,5\}$. As $K$ ranges over degree $d$ extensions of $\Q$ with discriminant $0<D_K \leq X$ (and non-$D_4$ in the case $d=4$),
\[  \sum_{\bstack{\deg(K)=d}{0<D_K \leq  X}}  |\Cl_K[\ell]|  \ll_{d,\ep} X^{\frac{3}{2} - \frac{1}{2\ell (d-1)}+ \ep}  ,\]
for all integers $\ell \geq \ell(d)$, where $\ell(2) = \ell(3)=1$, $\ell(4)=8$, $\ell(5)=25$.
\end{cor}

\begin{cor}\label{cor_main_small}
For positive integers $\ell \leq 7 $, averaging over non-$D_4$-quartic fields,
\[  \sideset{}{'}\sum_{\bstack{\deg(K)=4}{0 < D_K \leq X}}  |\Cl_K[\ell]|  \ll_{d,\ep} X^{\frac{3}{2}  - \frac{1}{48}+ \ep}  .\]
For positive integers $\ell \leq 24$, averaging over quintic fields,
\[   \sum_{\bstack{\deg(K)=5}{0 < D_K \leq X}}   |\Cl_K[\ell]|  \ll_{d,\ep} X^{\frac{3}{2}  - \frac{1}{200}+ \ep}  .\]
\end{cor}

Our strategy is as follows. Recall that  $N_d(X)$ denotes the number of degree $d$ fields $K$ over $\Q$, up to isomorphism, with $0 < D_K \leq X$, and let $N_d(X;p)$ denote the number of degree $d$ fields $K$ over $\Q$, up to isomorphism, with $0 <D_K \leq X$,  such that the rational prime $p$ splits completely in $K$. (For $d=4$ we define $\tilde{N}_4(X;p)$ analogously, restricting to non-$D_4$-quartic fields.)
Suppose we know that for each fixed prime $p$, $N_d(X;p)$ is a positive proportion of $N_d(X)$, so $p$ splits completely in a positive proportion of the fields. Then one would expect the fields in which the primes split completely to distribute somewhat evenly, so that ``most fields'' have the property that ``near the average number'' of primes split completely in them; that is, one would expect that the primes do not conspire to cause many fields to fail the criterion of Theorem \ref{lemma_M_primes}. We will make this argument precise by developing a flexible ``Chebyshev sieve'' (Lemma \ref{lemma_sieve}, related to Chebyshev's inequality); the crucial input to the sieve will be asymptotics for  $N_d(X;p)$ with power-saving error and explicitly given dependence on $p$ (Lemma \ref{lemma_field_2}, Theorem \ref{thm_field_3}, Theorems \ref{thm_field_4} and \ref{thm_field_5}). 

Counting quadratic fields may be accomplished by a simple classical argument (given in an appendix, Section \ref{sec_real_quad}).
Power-saving error terms for $N_d(X)$ were first found in the cases $d=3,4$ by Belabas, Bhargava, and Pomerance \cite{Bel10}, and first found in the case $d=5$ by Shankar and Tsimerman \cite{ShaTsi14}.
In the case $d=3,$ Bhargava, Shankar, and Tsimerman \cite{BST}
 and Taniguchi and Thorne \cite{TanTho13}  have also proved a second main term and improved the power-saving error term.
For the refined estimates that we require on $N_d(X;p)$, we quote the necessary asymptotics for $d=3$ from \cite{TanTho13}, while for $d=4,5$ we prove the necessary estimates using the methods and results from \cite{Bel10} and \cite{ShaTsi14}.
In fact, in Sections~\ref{sec_field_4} and \ref{sec_field_5}, we give the field counting asymptotics for fields with any chosen splitting types at a finite set of primes with the expectation that they could be useful in other applications; see Theorems \ref{thm_field_4_plus} and \ref{thm_field_5_plus}.
Our counting theorems improve upon analogous results that appear in
four recent papers, three \cite{Yan09,  ChoKim15, SST15}  in the area of finding symmetry groups of families 
of $L$-functions (see \cite{SarShinTem15} for a general overview of the area) and one 
\cite{LOTho14} studying the distribution of ramified primes in small-degree number fields.  See  Sections~\ref{sec_field_4} and \ref{sec_field_5} for detailed comparisons to these previous works.

\section{Anatomy of the proof}\label{sec_anatomy}

\subsection{Reduction to counting bad fields}\label{sec_anatomy_badfields}
We now outline the strategy in more detail; for the sake of motivation, we focus temporarily on proving upper bounds on average.
Let us fix $d$ and define for any degree $d$ field $K$ over $\Q$ and any real parameter $Y \geq 1$,
\[ N(K;Y)= \#\{ \text{rational primes $p \leq Y$ that split completely in $K$}\}. \]
(Implicitly,  in the case $d=4$ we further restrict to non-$D_4$-quartic fields.)
Let us fix a positive integer $\ell$ and a parameter $\del_1  < \frac{1}{2\ell (d-1)}$, to be chosen precisely later. Then by Theorem \ref{lemma_M_primes}, for any $X \geq 1$, 
\begin{eqnarray*}
\sum_{X <  D_K \leq 2X}  |\Cl_K[\ell]| 
	& \ll & \sum_{X <  D_K \leq 2X} D_K^{1/2 + \ep} N(K;D_K^{\del_1})^{-1} \\
	& \ll & X^{1/2 + \ep} \sum_{X <  D_K \leq 2X}  N(K;X^{\del_1})^{-1}.
	\end{eqnarray*}
Now given real parameters $X \geq 1$ and $1 \leq M \leq Y$, we define $\Bscr_d^0(X;Y,M)$ to be the set 
\[\Bscr^0_d(X;Y,M) = \{ K/\Q, \; \deg(K)=d, \;  X <D_K \leq 2X : \\ \text{ at most $M$ primes $p \leq Y$ split completely in $K$}\},
\]
 (with the usual further restriction in the case $d=4$).

We denote by $\pi(Y)$ the number of rational primes $p \leq Y$, and let us regard $1 \leq M \leq \pi(X^{\del_1})$ as fixed for the moment, to be specified later.
Then we may make the decomposition
\[ \sum_{X <  D_K \leq 2X}  |\Cl_K[\ell]|  
\ll  X^{1/2 + \ep} \left(\sum_{\bstack{X <  D_K \leq 2X}{K \not\in \Bscr^0_d(X;X^{\del_1}, M)}} N(K;X^{\del_1})^{-1}
	+\sum_{K \in \Bscr^0_d(X;X^{\del_1}, M)}  N(K;X^{\del_1})^{-1} \right).\]
Since $N(K;X^{\del_1}) \geq M$ if $K \not\in \Bscr^0_d(X;X^{\del_1}, M)$, we have
\[  \sum_{X <  D_K \leq 2X}  |\Cl_K[\ell]|   \ll
X^{1/2 + \ep} \left(\sum_{\bstack{X <  D_K \leq 2X}{K \not\in \Bscr^0_d(X;X^{\del_1}, M)}}  M^{-1}
	+ \sum_{K \in \Bscr^0_d(X;X^{\del_1}, M)} 1 \right),\]
and we may conclude that
\beq\label{CL_bound_0}
 \sum_{X <  D_K \leq 2X}  |\Cl_K[\ell]|  
 \ll X^{3/2+\ep}M^{-1} + \#  \Bscr^0_d(X;X^{\del_1}, M) X^{1/2 + \ep}. 
 \eeq
 Then upon defining (with the usual further restriction in the case $d=4$) the set
 \beq\label{Bset_dfn'}
 \Bscr_d(X;Y,M) = \{ K/\Q, \; \deg(K)=d, \;  0 <D_K \leq X :  \text{ at most $M$ primes $p \leq Y$ split completely in $K$}\},
 \eeq
we may trivially replace the expression $ \#  \Bscr^0_d(X;X^{\del_1}, M)$ in (\ref{CL_bound_0}) by $ \#  \Bscr_d(2X;X^{\del_1}, M)$ and only increase the right-hand side. 
 
We now suppose that we can bound from above the cardinality of the  ``bad set'' $\Bscr_d(2X;X^{\del_1}, M)$ for appropriate $\del_1$ and $M$. Note that one expects via the Chebotarev density theorem that a positive proportion of the primes up to $X^{\del_1}$ split completely in $K$, so that a reasonable choice for $M$ will be proportional to $\pi(X^{\del_1})$. Precisely, we suppose that there is a small fixed $\del_2>0$ such that for every $X\geq 1$ and an appropriate choice of $M$ with  $X^{\del_1}/\log X \ll M \ll X^{\del_1}/\log X$ we have
\beq\label{B_bound_explain}
\#  \Bscr_d(2X;X^{\del_1}, M) \ll X^{1-\del_2 + \ep},
\eeq
for all $\ep>0$.
 Then upon summing over $O(\log X)$ ranges and applying (\ref{CL_bound_0}) and (\ref{B_bound_explain}) within each range, we see that for any $X \geq 1$,
  \begin{eqnarray}\label{CL_bound}
 \sum_{0 <  D_K \leq X}  |\Cl_K[\ell]|  
& \leq &\sum_{0 \leq j \leq \lceil \log_2 X \rceil} \sum_{2^{j-1} < D_K \leq 2^{j}} |\Cl_K[\ell]| \nonumber \\
& \ll  & \sum_{0 \leq j \leq \lceil \log_2 X \rceil}  \left\{ (2^{j-1})^{3/2+\ep}(2^{(j-1) \del_1})^{-1}\log 2^j + \#  \Bscr_d(2^j;2^{(j-1)\del_1}, M) (2^{j-1})^{1/2 + \ep} \right\} \nonumber \\ 
& \ll  & \log X \sum_{0 \leq j \leq  \lceil \log_2 X \rceil}  \left\{ (2^j)^{3/2 - \del_1 +\ep} + (2^j)^{3/2 - \del_2 + 2\ep} \right\} \ll  X^{3/2 - \del +3\ep},
 \end{eqnarray}
 where $\del = \min\{\del_1,\del_2\}$ and $\ep>0$ is arbitrarily small. 
Thus we see that an upper bound of the form (\ref{B_bound_explain}) is the key to obtaining an average result in the shape of Corollaries \ref{cor_main_d} and \ref{cor_main_small}; this upper bound plays a similarly crucial role in obtaining the results of Theorem \ref{thm_main}, as we show in Section \ref{sec_main_theorems}.

Ultimately, we will prove the following version of (\ref{B_bound_explain}), which controls the number of possible bad fields:
\begin{thm}\label{thm_B}
Let  $ \Bscr_d(X;Y,M)$ be defined as in  (\ref{Bset_dfn'}).
Set 
\beq\label{del_choice_d_specific}
 \del_0(d) = \begin{cases}
	1/6 & \text{if $d=2$} \\
	2/25 & \text{if $d=3$}\\
	1/48 & \text{if $d=4$} \\
	1/200 & \text{if $d=5$}.
		\end{cases}
		\eeq
For each $d=2,3,4,5$, there is a constant $0<c_0(d)<1$ such that for every $0<\del \leq \del_0(d)$ and every $X \geq 1$,
\[ \# \Bscr_d(X;(X/2)^{\del}, \frac{1}{2}c_0(d)  \frac{(X/2)^{\del}}{\log (X/2)^{\del}}) \ll X^{1  - \del +\ep} \]
for every $\ep>0$.
\end{thm}
\begin{remark}
The methods of this paper also prove an analogous theorem if the condition ``split completely'' in the definition (\ref{Bset_dfn'}) is replaced by another fixed splitting type.
\end{remark}

\subsection{Counting bad fields via a sieve and counts for fields with local conditions}
We prove Theorem \ref{thm_B} via a sieve we develop for this purpose; to describe the strategy,
we first recall the simplest classical setting of a sieve. Let $\mathscr{A}$ be a finite set of elements of cardinality $N$, and let $\Pscr$ denote the set of all rational primes. We assume a certain property of interest has been specified so that each element $a \in \Ascr$ either satisfies it or not, with respect to $p$, for each $p \in \Pscr$.  For each prime $p \in \Pscr$ we let $\Ascr_p$ denote the finite subset of $\Ascr$ that satisfies the fixed property with respect to the prime $p$. Moreover we assume we know that for each $p$ there exists a real number $0 \leq \del_p<1$ and a real number $R_p$ with $|R_p| \leq N$ such that
\beq\label{A_inf} \# \Ascr_p = \del_p N + R_p. 
\eeq
In simplest terms, a classical aim of a sieve  is to provide an upper bound for the number of elements in the set $\Ascr$ such that the designated property fails for all primes $p \leq z$, for some fixed threshold $z$. Thus one could use a sieve to provide an upper bound for 
\[ \# ( \Ascr \setminus \Union_{p \leq z} \Ascr_p ).\]
For example, to sieve for prime numbers, the set $\mathscr{A}$ is a finite set of integers, and the property is that $p|a$.
Slightly more generally, one could apply a classical sieve such as  the Tur\'{a}n sieve to count 
\beq\label{Turan}
  \# ( \Ascr \setminus \Union_{p \in P_0} \Ascr_p )
  \eeq
for an arbitrary fixed finite set of primes $P_0$.

In our application, the set $\mathscr{A}$ is the set of fields $K/\Q$ of degree $d$ with $D_K \in (0,X]$ and the property is that $p$ splits completely in $K$, so that $\Ascr_p$ is the subset of fields in which the prime $p$ splits completely. 
 In this setting, assuming we possess an appropriate understanding of $\# \Ascr_p$ as in  (\ref{A_inf}), then (\ref{Turan}) would allow us to count those degree $d$ fields $K$ with $D_K \in (0,X]$ in  which a \emph{fixed} set of primes fail to split completely. But in order to bound the bad set $\Bscr_d(X;X^{\del_1}, M)$ we require more flexibility: a field belongs to this set if all the primes in a sufficiently large set fail to split completely in $K$, but the relevant large set of primes might be different for two different bad fields $K$. Thus we develop in Section \ref{sec_sieve} a flexible new sieve that allows us to count elements $a\in \Ascr$ that fail to lie in $\Ascr_p$ for many $p$, without specifying which $p$ fail for any given $a$.

The key input to any sieve is an understanding of $\Ascr_p$ that provides the expression (\ref{A_inf}). In our case, this requires an understanding of $N_d(X)$, $N_d(X;p)$, and $N_d(X;pq)$ for two distinct primes $p,q$; here $N_d(X;pq)$ counts the number of degree $d$ fields $K/ \Q$ in which both $p$ and $q$ split completely. In the case of quartic fields, we let $\tilde{N}_4(X)$, $\tilde{N}_4(X;p)$ and $\tilde{N}_4(X;pq)$ denote the analogous quantities, restricted to non-$D_4$-quartic fields $K/\Q$.

We now summarize the key results we will require for the sieve.
For quadratic fields, we record:
\begin{lemma}\label{lemma_field_2}
There exists a constant $c_2>0$, such that for $e=e_1$ or $e=e_1e_2$ for distinct primes $e_1,e_2$,
\begin{eqnarray}
N_2(X) & = & c_2 X + O(X^{1/2})  \label{N2}\\
N_2(X;e) & = & \del_{e} c_2 X + O(eX^{1/2}) \label{N2p}
\end{eqnarray}
where $\del_e$ is a multiplicative function defined for any prime $e$ by 
\beq\label{del_p_2}
\del_e = \frac{1}{2} \frac{1}{(1 + e^{-1})}.
\eeq
\end{lemma}
For completeness, we record a simple proof of this classical result in an appendix (Section \ref{sec_real_quad}); the error terms given here can be improved (see for example the survey \cite{Pap05}) but will suffice for our application.

In contrast, the results for cubic, quartic, and quintic fields are deep. 
For cubic fields, we  cite work of  Taniguchi and Thorne \cite{TanTho13}:
\begin{letterthm}[Theorems 1.1, 1.3 \cite{TanTho13}]\label{thm_field_3}
There exist constants $c_3>0$, $c_3'<0$ such that for $e=e_1$ or $e=e_1e_2$ for distinct primes $e_1,e_2$,
\begin{eqnarray}
N_3(X)  &=& c_3X + c_3'X^{5/6} + O(X^{7/9+\ep}) \label{N3}\\
N_3(X;e) & = & \del_e c_3X + \del_e'c_3'X^{5/6} + O(e^{8/9}X^{7/9+\ep}) \label{N3p}
\end{eqnarray}
where $\del_e$ and $\del'_e$ are multiplicative functions defined  for any prime $e$ by 
\beq\label{del_p_3}
\del_e = \frac{1}{6}\frac{1}{(1+e^{-1} + e^{-2})}, \qquad \del_e' = \frac{1}{6} + O(e^{-1/3}).
\eeq
\end{letterthm}

For quartic fields, we have:
\begin{thm}\label{thm_field_4}
There exists a constant $c_4>0$ such that for $e = e_1$ or $e=e_1e_2$ for distinct primes $e_1,e_2$,
\begin{eqnarray}
\tilde{N}_4(X) & = & c_4X + O(X^{23/24 + \ep}) \label{N4} \\
\tilde{N}_4(X;e) & = & \del_e c_4X  + O(e^{1/2 + \ep} X^{23/24+\ep}) \label{N4p}
\end{eqnarray}
where $\del_{e}$ is a multiplicative function defined for any prime $e$ by
\beq\label{del_p_4}
 \del_e = \frac{1}{24}\frac{1}{(1+ e^{-1} + 2e^{-2} + e^{-3})}.
 \eeq
\end{thm}
We note (\ref{N4}) is due to \cite[Theorem 1.3]{Bel10}; we deduce (\ref{N4p})  in Section \ref{sec_field_4}, using the methods of Belabas, Bhargava, and Pomerance \cite{Bel10}, which build on the work of Bhargava \cite{Bha05} that obtained the original count of $S_4$-quartic fields with an $o(X)$ error term. See Theorem \ref{thm_field_4_plus} for our most general result of this type, of which Theorem \ref{thm_field_4} is a special case.

For quintic fields, we have:

\begin{thm}\label{thm_field_5}
There exists a constant $c_5>0$ such that for $e = e_1$ or $e=e_1e_2$ for distinct primes $e_1,e_2$,
\begin{eqnarray}
N_5(X) & = & c_5X + O(X^{199/200 + \ep}) \label{N5}\\
N_5(X;e) & = & \del_e c_5X  + O(e^{1/2+\ep}X^{79/80+\ep}+X^{199/200+\ep})  \label{N5p}
\end{eqnarray}
where $\del_{e}$ is a multiplicative function defined for any prime $e$ by
\beq\label{del_p_5}
 \del_e = \frac{1}{120} \frac{1}{(1+ e^{-1} + 2e^{-2} + 2e^{-3}+e^{-4})}.
 \eeq
\end{thm}
We note (\ref{N5}) is due to \cite{ShaTsi14}; we deduce (\ref{N5p})  in Section \ref{sec_field_5}, using the methods of Shankar and Tsimerman \cite{ShaTsi14}, which  build on the work  of Bhargava  \cite{Bha10a} that obtained the original count of $S_5$-quintic fields  with an $o(X)$ error term.  (We also fill in a missing step from \cite{ShaTsi14}.)
See Theorem \ref{thm_field_5_plus} for our most general result, of which Theorem \ref{thm_field_5} is a special case.

We remark that the techniques for counting number fields that produced these results for $N_d(X;e)$ continue to be refined, and we may expect that the error terms will continue to be reduced. Thus in our subsequent computations involving $N_d(X;e)$ we have worked more generally with error terms of the form $O(e^\sig X^\tau)$, so that it will be immediately clear how improvements in counting fields will lead to refinements of our results.  (In particular,  improved error terms for \emph{smoothed} versions of the counting functions $N_d(X;e)$ would suffice for our application.) We note that the mechanism we employ will apply equally well to higher degree extensions of $\Q$ (or extensions of a fixed number field, using the more general form of Theorem \ref{lemma_M_primes} available in \cite[Lemma 2.3]{EllVen07})  if suitable results for $N_d(X)$ and $N_d(X;e)$ (or their analogues for extensions of a fixed number field) become available. In addition, one might consider other families of fields for which precise asymptotics are known, such as abelian fields over $\Q$ with a fixed  Galois group, ordered either by discriminant \cite{Wri89,FLN15} or by conductor \cite{Woo10}. It would be an interesting question to see whether the existing methods can be refined to produce an appropriate power-saving error term with sufficiently explicit dependence on a finite number of local conditions.

\section{The Chebyshev sieve}\label{sec_sieve}
We now develop in a fully general setting a new sieve that allows us to give an upper bound for the number of elements $a$ belonging to a set $\Ascr$ that satisfy a desired property with respect to $p$ for ``few'' $p$ (without specifying for which $p$ it is satisfied). We  will see that the principal idea is probabilistic, relating to Chebyshev's inequality, thus we dub it the Chebyshev sieve.

As before, let  $\mathscr{A}$ be a finite set of cardinality $N$,  let $\Pscr$ denote the set of all rational primes, and let $\Ascr_p$ denote the finite subset of $\Ascr$ that satisfies the fixed property with respect to the prime $p$. 
For a fixed real parameter $z \geq 1$, we let 
\[P(z)  = \prod_{\bstack{p \in \Pscr}{p \leq z}} p\]
and we define for each $a \in \Ascr$  the quantity
\[ N(a) = \# \{ p| P(z) : a \in \Ascr_p \}.\]
Next, we set 
\beq\label{M_dfn}
 M(z) = \frac{1}{N} \sum_{a \in \Ascr} N(a)= \frac{1}{N} \sum_{p|P(z)} \#\Ascr_p 
 \eeq
to be the mean  number of sets $\Ascr_p$ (with $p \leq z$) to which a typical element $a \in \Ascr$ belongs. (In non-vacuous cases,  $M(z)$ is nonzero.)
We would expect that a typical element $a \in \Ascr$ has $N(a)$ being about size $M(z)$, and we want to bound from above the number of $a \in \Ascr$ which have $N(a)$ being unusually small, that is, less than a fixed small proportion of $M(z)$.

Given $1 \leq M \leq z$, we define $\Escr(\Ascr; z,M)$ to be the set of elements $a \in \Ascr$ such that at most $M$ primes $p | P(z)$ have $a \in \Ascr_p$. (Or in other words, $\Escr(\Ascr; z,M)$ is the set of elements $a \in \Ascr$ such that $N(a) \leq M$.) Then we set 
\[E(\Ascr;z,M) = \# \Escr(\Ascr; z,M).\]
 Our sieve lemma will provide us with an upper bound for $E(\Ascr;z,\frac{1}{2}M(z))$; that is, the number of elements in $\Ascr$ that lie in $\Ascr_p$ for fewer than half the mean number of $p$.

For the purposes of the  lemma, we introduce the following notation.
Given distinct primes $p,q$  we let $\Ascr_{pq} = \Ascr_p \intersect \Ascr_q$, and let $R_{p,q}$ denote the quantity such that 
\[
 \# \Ascr_{pq} = \del_p \del_q N + R_{p,q}.
\]
(For notational convenience,  we will interpret $R_{p,p}$ as $R_p$.)
Finally, we set
\[ U(z) = \sum_{p | P(z)} \del_p.\]

We now state the key sieve lemma.
\begin{lemma}[Chebyshev Sieve]\label{lemma_sieve}
With the setting described above,
\[ E(\Ascr; z, \frac{1}{2}M(z)) \leq 
  \frac{4N}{M(z)^2} \left( U(z)+  \frac{1}{N}  \sum_{p,q| P(z)} |R_{p,q}|  
  +  \frac{2U(z)}{N} \sum_{p | P(z)} |R_p|
   + \left( \frac{1}{N} \sum_{p | P(z)} |R_p| \right)^2 \right).
\]
\end{lemma}

\subsection{Proof of the sieve lemma}\label{sec_sieve_proof}
We note that the sieve inequality we prove is related to the classical Tur\'{a}n sieve (see for example Theorem 4.1.1 of \cite{CojMur06}), and can be seen as an application of Chebyshev's inequality
\[\P(|X-\mu|\geq \al)\leq \sigma^2/\al^2,\]
 for $X$ a random variable with mean $\mu$ and variance $\sigma^2$, applied to the random variable $N(a)$ when $a$ is drawn uniformly from $\Ascr$.
 
We prove the lemma directly.
We begin by noting that
\[ \frac{1}{N} E(\Ascr; z, \frac{1}{2}M(z)) (\frac{1}{2} M(z))^2
 	  \leq   \frac{1}{N} \sum_{a \in \Escr(\Ascr; z, \frac{1}{2}M(z))}  (N(a) - M(z))^2 
	 \leq  \frac{1}{N} \sum_{a \in \Ascr} (N(a) - M(z))^2.  
\]
It then suffices to prove the variance term on the right-hand side satisfies
\beq \label{variance}
\frac{1}{N} \sum_{a \in \Ascr} (N(a) - M(z))^2
	\leq  U(z)+ \frac{1}{N}  \sum_{p,q| P(z)} |R_{p,q}|  + 2U(z)\left( \frac{1}{N} \sum_{p | P(z)} |R_p| \right)  + \left( \frac{1}{N} \sum_{p | P(z)} |R_p| \right)^2.
\eeq	
We first note from (\ref{M_dfn}) that the mean satisfies
\beq\label{mean}
M(z)= \frac{1}{N} \sum_{p|P(z)} \#\Ascr_p
	= \frac{1}{N} \sum_{p| P(z)} (\del_p N + R_p)
	= U(z) + \frac{1}{N} \sum_{p | P(z)} R_p.
\eeq
We now consider the left-hand side of (\ref{variance}), which we trivially expand as
\beq\label{expand}
\frac{1}{N} \sum_{a \in \Ascr} N(a)^2 -  \frac{2}{N}\sum_{a \in \Ascr} N(a) M(z)
		+ M(z)^2 = \frac{1}{N} \sum_{a \in \Ascr} N(a)^2 -   M(z)^2.
		\eeq
The first term on the right-hand side of (\ref{expand}) is equal to
\begin{eqnarray*}
 \frac{1}{N} \sum_{p,q | P(z)} \#( \Ascr_p \intersect \Ascr_q) 
	& = & \frac{1}{N}  \left( \sum_{p |P(z)} \del_p N +  \sum_{\bstack{p,q| P(z)}{p \neq q}} \del_p \del_q N +  \sum_{p,q| P(z)} R_{p,q} \right) \\
	& = & \sum_{p | P(z)} \del_p + \left( \sum_{p | P(z)} \del_p \right)^2 - \sum_{p | P(z)} \del_p^2 + \frac{1}{N}  \sum_{p,q| P(z)} R_{p,q}  \\
	& = &  \sum_{p | P(z)} \del_p(1-\del_p) + U(z)^2 + \frac{1}{N}  \sum_{p,q| P(z)} R_{p,q}.
	\end{eqnarray*}
On the other hand, we may expand $M(z)^2$ via (\ref{mean}) and see that after cancellation of the $U(z)^2$ factor, the right-hand side of (\ref{expand}) is equal to
\[ \sum_{p | P(z)} \del_p(1-\del_p)  + \frac{1}{N}  \sum_{p,q| P(z)} R_{p,q}  - 2U(z)\left( \frac{1}{N} \sum_{p | P(z)} R_p \right)  - \left( \frac{1}{N} \sum_{p | P(z)} R_p \right)^2.\]
As $R_p$ may be either positive or negative, we take absolute values; then using the fact that $\del_p\leq 1$ we see the resulting inequality simplifies to  (\ref{variance}), thus proving the lemma.

\section{Asymptotic count of non-$D_4$-quartic fields}\label{sec_field_4}

In this section we will prove the following, of which Theorem~\ref{thm_field_4} is a special case.
\begin{thm}\label{thm_field_4_plus}
Let $P$ be a finite set of primes.  For each prime $p\in P$ we choose a splitting type at $p$ and assign a corresponding density as follows:
\begin{align*}
\del_p&:= \frac{1}{24}(1+ p^{-1} + 2p^{-2} + p^{-3})^{-1} &\textrm{ for }p=\wp_1\wp_2\wp_3\wp_4\\
\del_p&:= \frac{1}{4}(1+ p^{-1} + 2p^{-2} + p^{-3})^{-1}  &\textrm{ for }p=\wp_1\wp_2\wp_3\\
\del_p&:= \frac{1}{3}(1+ p^{-1} + 2p^{-2} + p^{-3})^{-1}  &\textrm{ for }p=\wp_1\wp_2 \textrm{ with $\wp_2$ inertia degree 3}\\
\del_p&:= \frac{1}{8}(1+ p^{-1} + 2p^{-2} + p^{-3})^{-1}  &\textrm{ for }p=\wp_1\wp_2 \textrm{ with $\wp_i$ inertia degree 2}\\
\del_p&:= \frac{1}{4}(1+ p^{-1} + 2p^{-2} + p^{-3})^{-1}  &\textrm{ for }p=\wp_1\\
\del_p&:= \frac{p^{-1} + 2p^{-2} + p^{-3}}{(1+ p^{-1} + 2p^{-2} + p^{-3})} &\textrm{ for $p$ ramified}.
\end{align*}
Let $\delta_P:=\prod_{p\in P} \delta_p$ and let $e=\prod_{p\in P} p$.
Let $\tilde{N}_4(X;P)$ be the number of non-$D_4$ quartic fields with absolute discriminant at most $X$ such that for each $p\in P$, the prime $p$ splits in the quartic field in the splitting type chosen for $p$ above.
There exists a constant $c_4>0$ such that
\begin{eqnarray}
\tilde{N}_4(X;P) & = & \del_P c_4X  + O(e^{1/2 + \ep} X^{23/24+\ep}), 
\end{eqnarray}
where the implied constant in the $O$ term is absolute (does not depend on $P$).
Moreover,  we may choose more than one splitting type at each prime and let $\delta_p$ be the sum of the corresponding densities and the result still holds.
\end{thm}
Bhargava \cite{Bha05} first determined the asymptotic count of non-$D_4$-quartic fields, and 
 Belabas, Bhargava, and Pomerance \cite{Bel10} gave a power-saving asymptotic for this count.
 We will follow the method of \cite{Bel10}, additionally requiring our chosen splitting types.  While the main term for such a restricted count appears in \cite[Theorem 3]{Bha05} (at least for one prime, and the same argument would work for more primes), we require a power-saving error term \emph{with explicit dependence} on the primes.   In fact, such results have appeared at least four times recently, but we will improve upon the exponents in all of these results and remove various hypotheses that don't hold in the situation in which we need to apply the bound.
 Yang \cite[Proposition 3.1.7]{Yan09} proved such a power-saving error of the form $\tilde{N}_4(X;P)  =  \del_P c_4X  + O(e^{2} X^{143/144+\ep})$. (\cite[Proposition 3.1.7]{Yan09} only states this for one local condition, but \cite[Section 7]{ChoKim15} remarked it can be extended to finitely many local conditions.)
Lemke Oliver and Thorne \cite[Theorem 2.1]{LOTho14} proved a power-saving error (in which we may only specify that $p$ is ramified)  of $\tilde{N}_4(X;P)  =  \del_P c_4X  + O(e^{9/10} X^{239/240+\ep}).$
Shankar, S\"{o}dergren, and Templier \cite{SST15} proved $\tilde{N}_4(X;P)  =  \del_P c_4X  + O(e^{12} X^{23/24+\ep})$ when $P$ contains a single prime.
 
  The exposition of the method in \cite{Bha05} and \cite{Bel10} is quite clear, so we will focus here on the particular aspects of the computation we need.    Instead of directly counting quartic fields, the method, equivalently, counts maximal quartic orders.
 The parametrization of quartic rings with their cubic resolvents due to Bhargava \cite{Bha04} (see also
\cite[Theorem 4.1]{Bel10}) gives an injection from the set of isomorphism classes of maximal quartic orders to the set of $\GL_2(\Z)\times\SL_3(\Z)$ classes of pairs of ternary quadratic forms with integral coefficients.  Pairs of integral ternary quadratic forms comprise a $12$ dimensional lattice $V_{\Z}=\Z^{12}$. Counting $\GL_2(\Z)\times\SL_3(\Z)$  classes of lattice points in $\Z^{12}$ is the same as
  counting lattice points in a fundamental domain for $\GL_2(\Z)\times\SL_3(\Z)$ on $\R^{12}$.
 In this paper, we need to count only these lattice points in particular translates of sublattices of $\Z^{12}$.
 We  collect some basic facts about the lattice translates corresponding to our desired fields, apply the geometry of numbers result from \cite{Bel10} to count the necessary lattice points, and then work to minimize the resulting error terms.

  As in \cite[Section 2.2]{Bha05} and \cite[Section 4]{Bel10} we use a certain random fundamental domain for the action of $\GL_2(\Z)\times\SL_3(\Z)$ on $\R^{12}$.  
For a positive integer $m$, let $L$ be a translate $v+mV_\Z$ ($v\in V_\Z$) of the sublattice $m V_\Z$ of $ V_\Z$. 
  Let $N'(L;X)$ denote the expected number of lattice points in $L$, with first coordinate non-zero and discriminant less than $X$, in a random fundamental domain.   (This notion of expected value for a random fundamental domain is defined as in  \cite[Equation (5)]{Bha05}, with $S$ the set of points of $L$ with first coordinate non-zero, but \emph{without} the ``abs. irr.'' condition that appears in \cite[Equation (5)]{Bha05}. See also \cite[p. 198]{Bel10}.) 
Let $N_{S_4}(q;X)$ be the number of classes in $V_\Z$ corresponding to isomorphism classes of $S_4$-quartic 
orders and whose index in their maximal order is divisible by $q$ and whose discriminant is less than $X$.
We have the following result that estimates these counts.

\begin{letterthm}[Theorem 4.11 of \cite{Bel10}]\label{T:countintrans}
Let $L$ be a translate $v+mV_\Z$ ($v\in V_\Z$).  Let $(a,b,c,d)$ denote the smallest positive first four coordinates of any element of $L$.
Then
$$
N'(L,X) = \frac{N_{S_4}(1;X)}{m^{12}} +O( \sum_S  \frac{X^{(|S|+\alpha_S+\beta_S+\gamma_S+\delta_S)/12}}{m^{|S|}a^{\alpha_S}b^{\beta_S}c^{\gamma_S}d^{\delta_S}} +\log X )
$$
where $S$ ranges over the non-empty proper subsets of the set of $12$ coordinates on $V_\Z$, and $\alpha_S,\beta_S,\gamma_S,\delta_S\in[0,1]$ are real constants that depend only on $S$ and satisfy $|S|+\alpha_S+\beta_S+\gamma_S+\delta_S\leq 11$.   
\end{letterthm}

Let $q$ be square-free and $(q,e)=1$.  First, we will assume that we have chosen unramified splitting types at each prime in $P$.
Now, we will start by counting the expected number $N'(q,e;X)$ of lattice points in a random fundamental domain that satisfy the following conditions: (1) their first coordinate is non-zero, (2) their discriminant is less than $X$, (3) their corresponding quartic ring is not maximal at each prime dividing $q$ and is maximal and of chosen splitting type at primes in $P$.
 We do this by summing Theorem~\ref{T:countintrans} over the collection $T$ of translates of $eq^2 V_\Z$ that give quartic rings that are not maximal at each prime dividing $q$, and are maximal and with chosen local splitting at each $p\in P$.
(See \cite[Section 4]{Bha04} for a description of which pairs of ternary quadratic forms correspond to quartic rings that are maximal or split in a certain way at a prime.)

Given $(a,b,c,d)\in [1,eq^2]^4$, we need to bound the number of translates in $T$ that have $(a,b,c,d)$ as the smallest positive first four coordinates of any element.  By \cite[Corollary 4.8]{Bel10}, there are $O(6^{\omega(q)}q^{14})$ translates of $q^2V_{\Z}$ that
are congruent to $(a,b,c,d)$ modulo $q^2$
and whose lattice points
 correspond to quartic rings that are not maximal at each prime dividing $q$.  Since $V_{\Z}$ is 12 dimensional, there are $e^8$ translates of $eV_{\Z}$ congruent to $(a,b,c,d)$ modulo $e.$  Thus by the Chinese Remainder Theorem, there are 
 $O(6^{\omega(q)}q^{14}e^8)$ 
translates in $T$ that have $(a,b,c,d)$ as the smallest positive first four coordinates of any element. 
 
For $q$ square-free, we define $\nu(q)$ to be the multiplicative function defined for a prime $p$ by
$$
\nu(p):=p^{-2}+2p^{-3}+2p^{-4}-3p^{-5}-4p^{-6}-p^{-7}+3p^{-8}+3p^{-9}-p^{-10}-p^{-11}.
$$
This is the density of lattice points that correspond to quartic rings non-maximal at $p$ \cite[Lemma 4.4]{Bel10}.
Then 
$\#T=\nu(q)q^{24} e^{12} \Gamma_P$, where
$$
\Gamma_P:=\prod_{p\in P} \delta_{p} (1-\nu(p)), 
$$
and $0 \leq \delta_{p} \leq 1$ is the density of lattice points corresponding to quartic rings that are
split as we chose at $p$ as a subset of those corresponding to quartic rings that are maximal at $p$ \cite[Lemma 23]{Bha04}.

If $q^2>X$, then all the classes counted by $N'(q,e;X)$ have discriminant $0$, and by  \cite[Lemma 4.10]{Bel10}, in this case there are $O(X^{11/12+\ep})$ such classes.

So now we consider the case when $q^2\leq X$, in which case by Theorem \ref{T:countintrans},
\begin{align*}
&N'(q,e;X)
\\
& \qquad =\nu(q) \Gamma_P {N_{S_4}(1;X)} 
+ O\left( \sum_{(a,b,c,d)\in [1,eq^2]^4}  6^{\omega(q)}q^{14}e^8 \left(  \sum_S  \frac{X^{(|S|+\alpha_S+\beta_S+\gamma_S+\delta_S)/12}}{(eq^2)^{|S|}a^{\alpha_S}b^{\beta_S}c^{\gamma_S}d^{\delta_S}} +\log X \right) \right).
\end{align*}
 
We have
\begin{align*}
&\sum_{(a,b,c,d)\in [1,eq^2]^4}  6^{\omega(q)}q^{14}e^8 \left(  \sum_S  \frac{X^{(|S|+\alpha_S+\beta_S+\gamma_S+\delta_S)/12}}{(eq^2)^{|S|}a^{\alpha_S}b^{\beta_S}c^{\gamma_S}d^{\delta_S}} +\log X \right) \\
&=  6^{\omega(q)}q^{14}e^8 \left( e^4q^8\log X + \sum_S  \frac{X^{(|S|+\alpha_S+\beta_S+\gamma_S+\delta_S)/12}}{(eq^2)^{|S|}} \sum_{(a,b,c,d)\in [1,eq^2]^4} \frac{1}{a^{\alpha_S}b^{\beta_S}c^{\gamma_S}d^{\delta_S}} \right) \\
&\leq  6^{\omega(q)}q^{14}e^8 \left(e^4q^8\log X+  \sum_S  \frac{X^{(|S|+\alpha_S+\beta_S+\gamma_S+\delta_S)/12}}{(eq^2)^{|S|}} 
\left( (eq^2)^{4-\alpha_S-\beta_S-\gamma_S-\delta_S} \log^4(eq^2)  \right)
  \right) \\
  &=  6^{\omega(q)}q^{22}e^{12} \left( \log X+ \sum_S  (X^{1/12}e^{-1}q^{-2})^{(|S|+\alpha_S+\beta_S+\gamma_S+\delta_S)}
\log^4(eq^2)  
 \right). \\
\end{align*}

Since $0\leq |S|+\alpha_S+\beta_S+\gamma_S+\delta_S \leq 11$, and recalling that $q^2 \leq X$, the above is 
\begin{align*}
&= O\left(  6^{\omega(q)}q^{22}e^{12}  \left( (X^{1/12}e^{-1}q^{-2})^{11} \log^4(eq^2) + \log^4(eq^2)  +\log X
\right)
 \right)\\
 &= O\left(  e^{1+\ep}  X^{11/12+\ep}  + q^{22}e^{12+\ep} X^{\ep} 
 \right).
\end{align*}

Let $N_{S_4}(q,e;X)$ be the number of classes in $V_\Z$, 
or equivalently lattice points in a fundamental domain,
corresponding to isomorphism classes of $S_4$-quartic orders, whose index in their maximal order is divisible by $q$ and whose discriminant is less than $X$, and that are maximal and of chosen splitting type at $p\in P$. 
Now, by inclusion-exclusion, as in the proof of \cite[Theorem 4.13]{Bel10}, we have that the number of isomorphism classes of maximal $S_4$-quartic orders splitting as chosen for $p\in P$ and having (absolute) discriminant less than $X$ is given by
$$
\sideset{}{'}\sum_{q\geq 1} \mu(q) N_{S_4}(q,e;X)
$$
where the sum is restricted to square-free $q$ that are relatively prime to $e$.

Now we compare $N_{S_4}(q,e;X)$ and $N'(q,e;X)$.   Note that the difference is that  $N'(q,e;X)$ excludes those lattice points with first coordinate $0$, and $N_{S_4}(q,e;X)$ excludes those lattice points that do not correspond to orders in $S_4$-quartic fields.  
So by \cite[Lemmas 4.9 and 4.10]{Bel10}, we have 
\[|N_{S_4}(q,e;X)-N'(q,e;X)| =O(X^{11/12+\ep}).\]
Thus by our previous computation for $N'(q,e;X)$,
\beq\label{sum_qN}
N_{S_4}(q,e;X)
=\nu(q) \Gamma_P {N_{S_4}(1;X)} +O\left(  e^{1+\ep}  X^{11/12+\ep}  + q^{22}e^{12+\ep} X^{\ep} 
 \right).
\eeq
So for a fixed $Q$ (to be chosen in terms of $X,e$ later), we sum over square-free $q$ with $(q,e)=1$ as in (\ref{sum_qN}), obtaining
\begin{align*}
\sideset{}{'}\sum_{q\geq 1} \mu(q) N_{S_4}(q,e;X) &= \sideset{}{'}\sum_{1\leq q\leq Q} \mu(q) N_{S_4}(q,e;X)  + \sideset{}{'}\sum_{q> Q} \mu(q) N_{S_4}(q,e;X) \\
&= \sideset{}{'}\sum_{1\leq q\leq Q} \mu(q) \nu(q) \Gamma_P {N_{S_4}(1;X)}
+O(E_1)
 + O(E_2) \\
 &= \sideset{}{'}\sum_{q \geq 1} \mu(q) \nu(q)\Gamma_P{N_{S_4}(1;X)}
+O(E_1)
 + O(E_2) +O(E_3)\\
 &= \prod_p (1-\nu(p)) \prod_{p\in P} \delta_{p} {N_{S_4}(1;X)}
+O(E_1)
 + O(E_2) +O(E_3)\\
\end{align*}

where
\begin{align*}
E_1&=  \sideset{}{'}\sum_{1\leq q\leq Q}  \left(  e^{1+\ep}  X^{11/12+\ep}  + q^{22}e^{12+\ep} X^{\ep} 
 \right)\\
E_2&= \sideset{}{'}\sum_{q> Q}\mu(q) N_{S_4}(q,e;X)\\
E_3&= \sideset{}{'}\sum_{q>Q} \nu(q) \Gamma_P {N_{S_4}(1;X)}.
\end{align*}
(Note that we handle the terms slightly differently than in \cite{Bel10}, so that $E_3$ above does not correspond to their $E_3$ term.)

We have $E_1=O( e^{1+\ep}  Q  X^{11/12+\ep} +Q^{23}e^{12+\ep}  X^{\ep}  )$.
By \cite[Lemma 4.3]{Bel10}, we have $N_{S_4}(q,e;X)=O(Xq^{-2+\ep})$, and so
$
E_2 =O(XQ^{-1+\ep}).  
$
We have $E_3=O(Q^{-1+\ep} X)$, since by \cite[Lemma 4.2]{Bel10}, we have $N_{S_4}(1;X)=O(X)$, and by definition $\nu(q) = O(q^{-2+\ep})$.  

If $e\leq X^{1/12}$, then we take $Q=X^{1/24}e^{-1/2}$, and we have 
\begin{align*}
\sideset{}{'}\sum_{q\geq 1} \mu(q) N_{S_4}(q,e;X) 
 &= \prod_p (1-\nu(p)) \prod_{p\in P} \delta_{p} {N_{S_4}(1;X)}
+O(e^{1/2+\ep}X^{23/24+\ep}).
\end{align*}
By \cite[Lemma 4.2]{Bel10}, we have that 
\[\prod_p (1-\nu(p))  {N_{S_4}(1;X)}=c_4X+O(X^{23/24+\ep}),\]
 for some positive constant $c_4$.
Thus we conclude that
the number of isomorphism classes of maximal $S_4$-quartic orders with our chosen splitting types at $p\in P$ and having (absolute) discriminant less than $X$ is 
$$\delta_Pc_4X 
+O(e^{1/2+\ep}X^{23/24+\ep}).$$

If $e>X^{1/12}$, then the number of isomorphism classes of maximal $S_4$-quartic orders with chosen splitting types for $p\in P$ and having (absolute) discriminant less than $X$ is $O(X)$ by \cite[Lemma 4.2]{Bel10}, which we may then also write as
$$\delta_Pc_4X 
+O(e^{1/2+\ep}X^{23/24+\ep}).$$

There are at most $O(X^{7/8 + \ep})$ quartic extensions with $D_K < X$ with Galois closure having Galois group $C_4, K_4$ or $A_4$ (Baily \cite{Bai80}, and Wong \cite{Won99b}).
So we can conclude Theorem~\ref{thm_field_4_plus} holds for unramified splitting types. This argument shows we can also choose more than one splitting type at each $p$, and sum the corresponding densities.

Now, given $P$ and choices for local splitting types some of which may be ramified, let $P_1$ be the subset of $P$ for which we choose only unramified splitting types. 
 We can find $\tilde{N}_4(X;P_1)$ using the result already proven.  For any subset $P_2\sub P\setminus P_1$, write
 $\tilde{N}_4(X;P_1 \cup \bar{P_2})$ for the number of non-$D_4$ quartic fields with absolute discriminant at most $X$ such that for each $p\in P_1$ the prime $p$ splits in one of our chosen spitting type, and for each $p\in P_2$ the prime $p$ \emph{does not} split in one of our chosen splitting types.  We can also apply the result already proven to find $\tilde{N}_4(X;P_1 \cup \bar{P_2})$.  Then using inclusion exclusion, we have
 \begin{align*}
 \tilde{N}_4(X;P)&=\sum_{P_2\sub P\setminus P_1} (-1)^{|P_2|} \tilde{N}_4(X;P_1 \cup \bar{P_2})\\
 &=\delta_Pc_4X + \sum_{P_2\sub P\setminus P_1} (-1)^{|P_2|} O(e^{1/2+\ep} X^{23/24+\ep}).
\end{align*}
Since each set $P_2$ corresponds to a distinct divisor of $e$ there are $O(e^{\ep})$ terms in the sum and Theorem~\ref{thm_field_4_plus} follows.

\begin{remark}\label{rem_d4} On the other hand, the number of $D_4$-quartic fields with $D_K<X$ is $\sim cX$ with $c\approx 0.052326...$, as initially indicated (as an order of magnitude) by Baily \cite{Bai80} and refined with an explicit constant by Cohen, Diaz y Diaz and Olivier \cite{CDO02}.  It is an interesting open problem to count $D_4$ fields with local conditions such as certain primes being split completely, and for now we exclude them from our consideration.
\end{remark}

\section{Asymptotic count of quintic fields}\label{sec_field_5}

In this section we will prove the following, of which Theorem~\ref{thm_field_5} is a special case.
\begin{thm}\label{thm_field_5_plus}
Let $P$ be a finite set of primes.  For each prime $p\in P$ we choose a splitting type at $p$ and assign a corresponding density as follows:
\begin{align*}
\del_p&:= \frac{1}{120}(1+ p^{-1} + 2p^{-2} + 2p^{-3}+p^{-4})^{-1} &\textrm{ for }p=\wp_1\wp_2\wp_3\wp_4\wp_5\\
\del_p&:=\frac{1}{12}(1+ p^{-1} + 2p^{-2} + 2p^{-3}+p^{-4})^{-1}&\textrm{ for }p=\wp_1\wp_2\wp_3\wp_4\\
\del_p&:= \frac{1}{8}(1+ p^{-1} + 2p^{-2} + 2p^{-3}+p^{-4})^{-1} &\textrm{ for }p=\wp_1\wp_2\wp_3 \textrm{ with $\wp_2,\wp_3$ inertia degree 2}\\
\del_p&:= \frac{1}{6}(1+ p^{-1} + 2p^{-2} + 2p^{-3}+p^{-4})^{-1} &\textrm{ for }p=\wp_1\wp_2\wp_3 \textrm{ with $\wp_3$ inertia degree 3}\\
\del_p&:= \frac{1}{6}(1+ p^{-1} + 2p^{-2} + 2p^{-3}+p^{-4})^{-1} &\textrm{ for }p=\wp_1\wp_2 \textrm{ with $\wp_2$ inertia degree 3}\\
\del_p&:= \frac{1}{4}(1+ p^{-1} + 2p^{-2} + 2p^{-3}+p^{-4})^{-1} &\textrm{ for }p=\wp_1\wp_2 \textrm{ with $\wp_2$ inertia degree 4}\\
\del_p&:= \frac{1}{5}(1+ p^{-1} + 2p^{-2} + 2p^{-3}+p^{-4})^{-1}  &\textrm{ for }p=\wp_1\\
\del_p&:= \frac{p^{-1} + 2p^{-2} + 2p^{-3}+p^{-4}}{(1+ p^{-1} + 2p^{-2} + 2p^{-3}+p^{-4})} &\textrm{ for $p$ ramified}.
\end{align*}
Let $\delta_P:=\prod_{p\in P} \delta_p$ and let $e=\prod_{p\in P} p$.
Let $N_5(X;P)$ be the number of quintic fields with absolute discriminant at most $X$ such that for each $p\in P$, the prime $p$ splits in the quartic field in the splitting type chosen for $p$ above.
There exists a constant $c_5>0$ such that
\begin{eqnarray}
N_5(X;P) & = & \del_P c_5X  + O(e^{1/2+\ep}X^{79/80+\ep}+X^{199/200+\ep}),
\end{eqnarray}
where the implied constant in the $O$ term is absolute (does not depend on $P$).
Moreover,  we may choose more than one splitting type at each prime and let $\delta_p$ be the sum of the corresponding densities and the result still holds.
\end{thm}

 Bhargava \cite{Bha10a} gave the first asymptotic count 
of quintic number fields, and Shankar and Tsimerman \cite{ShaTsi14}, building on Bhargava's work, gave the first power-saving error term.   
Both proofs fundamentally rely on Bhargava's \cite{Bha08a} parametrization of quintic rings.  
We will follow the outline of the argument of \cite{ShaTsi14}, additionally requiring our chosen splitting conditions.  While the main term for such a restricted count appears in \cite[Theorem 3]{Bha10a} (at least for one prime, and the same argument would work for more primes), we require a power-saving error term \emph{with explicit dependence} on the primes.   While such bounds have appeared in at least three recent papers, we will improve on the exponents in all of them, as well as remove hypotheses that do not hold in our cases of interest.
Lemke Oliver and Thorne \cite[Theorem 2.1]{LOTho14} have shown, assuming that we choose ramification at each prime $p\in P$, that $N_5(X;P)  =  \del_P c_5X  + O(eX^{199/200+\ep}).$
Cho and Kim \cite[Section 6]{ChoKim15} have recently proven a bound of the sort we desire; it seems they show
 $N_5(X;P)  =  \del_P c_5X  + O(e^{2-\ep}X^{399/400+\ep}) $.
Also, 
Shankar, S\"{o}dergren, and Templier \cite{SST15} stated the bound $N_5(X;P)  =  \del_P c_5X  + O(e^{40}X^{79/80+\ep}+X^{199/200+\ep})$ when $P$ contains a single prime.

 Instead of directly counting quintic fields, the method, equivalently, counts maximal quintic orders.
Analogously to the quartic case, we use a parametrization of quintic rings with their sextic resolvents due to Bhargava \cite{Bha08a}.  Let $V_\Z=\Z^{40}$ denote the space of quadruples of $5\times 5$ skew-symmetric matrices with integer coefficients.  Then quintic rings with their sextic resolvents are parametrized by 
$\GL_4(\Z)\times\SL_5(\Z)$ orbits on $V_\Z$ \cite[Theorem 1]{Bha08a}.  These orbits correspond to lattice points in a fundamental domain for
 $\GL_4(\Z)\times\SL_5(\Z)$ on $\R^{40}$.
   As in \cite[Section 2.2]{Bha10a} and \cite[Section 2.2]{ShaTsi14}, we take a certain random fundamental domain for the action of $\GL_4(\Z)\times\SL_5(\Z)$ on $\R^{40}$.
For a subset $S\sub V_\Z$, let $N_{\operatorname{dom}}(S;X)$ denote the expected number of elements of $S$ with absolute discriminant less than $X$ and whose associated quintic ring is an integral domain (i.e. is an order in a quintic field), in a random fundamental domain (as in \cite[Equation (1)]{ShaTsi14}, summed over the implicit $i$ there).  Let $N^*(S;X)$ denote the expected number of elements of $S$ with absolute discriminant less than $X$, in a random fundamental domain (as in the equation after (1) in \cite{ShaTsi14}, summed over the implicit $i$ there).

We first consider the case in which only unramified splitting types are chosen.
Let $a_{12}$ denote the $(1,2)$ coordinate of the first matrix in a quadruple of $5\times 5$ skew-symmetric matrices.
For a square-free integer $q$ relatively prime to $e$, let $W_{q,e}\sub V_\Z$ denote the set of elements corresponding to quintic rings that are not maximal at each prime dividing $q$ and are maximal and  of chosen splitting type at primes dividing $e$.  Recall from  \cite[Section 12]{Bha08a} that $W_{q,e}$ is defined by congruence conditions modulo $q^2e$ (for maximality, an argument analogous to that in \cite[Lemma 22]{Bha04} is necessary).
Let $U_e\sub V_\Z$ denote the set of elements corresponding to quintic rings that are  maximal at all primes and  of chosen splitting type  at the primes dividing $e$. Then counting $N_{\mathrm{dom}}(U_e;X)$ will provide us with precisely the count $N_5(X;e)$ we require.
We will count lattice points in $U_e$ by using inclusion-exclusion to reduce to counting lattice points in the $W_{q,e}$.

By \cite[Equation (27)]{Bha10a} (see also \cite[Equation (4)]{ShaTsi14}), if $L$ is a translate of the lattice $mV_\Z$ and $m=O(X^{1/40})$, then
\begin{equation}\label{E:quintictranslate}
N^*(L\cap \{a_{12}\ne 0 \};X)=c_0m^{-40}X+O(m^{-39} X^{39/40}),
\end{equation}
for some positive absolute constant $c_0$.

Bhargava gives the density of lattice points corresponding to rings maximal at a given prime \cite[Equation (48)]{Bha08a} and
the density of lattice points corresponding to rings maximal and of each splitting type \cite[Lemma 20]{Bha08a}.
Using these two computed densities, we conclude that  
of the $(q^2e)^{40}$ quadruples of $5\times 5$ skew-symmetric matrices  mod $q^2e$, we have that $W_{q,e}$
 corresponds to 
$
\nu(q)q^{80} \delta_P e^{40}\prod_{p\in P} (1-\nu(p))
$ of them,
where $$
\nu(p)=1-\frac{(p-1)^8p^{12}(p+1)^4(p^2+1)^2(p^2+p+1)^2(p^4+p^3+p^2+p+1)(p^4+p^3+2p^2+2p+1)}{p^{40}},
$$
and we extend this to a multiplicative function $\nu(q)$ for square-free $q$. (Here, $\nu(p)$ is the density of lattice points correspond to rings that are  \emph{non-maximal} at $p$ from \cite[Equation (48)]{Bha08a}.) 
Note that  $\nu(p)=p^{-2}+O(p^{-3}) $ and thus $\nu(q)= O(q^{-2+\ep})$. 

  We have that $\delta_P\leq 1$ and $1-\nu(p)\leq 1$.  
So,  when $q^2e=O(X^{1/40})$,  by summing Equation~\eqref{E:quintictranslate} over all the translates of $q^2eV_\Z$ that comprise $W_{q,e}$, we find that
\begin{align}\label{e:Nstar}
&N^*(W_{q,e}\cap \{a_{12}\ne 0 \};X)\\&=\nu(q)q^{80} \delta_P e^{40}\prod_{p\in P} (1-\nu(P))c_0q^{-80}e^{-40}X+O(\nu(q)q^{80} \delta_P e^{40}\prod_{p\in P} (1-\nu(p)) q^{-78}e^{-39} X^{39/40})\notag\\
&=\nu(q) \delta_P \prod_{p\in P} (1-\nu(p))c_0X+O(\nu(q)q^{2} \delta_P e\prod_{p\in P} (1-\nu(p))  X^{39/40})\notag\\
&=\nu(q) \delta_P \prod_{p\in P} (1-\nu(p))c_0X+O(q^{\ep}  e  X^{39/40}),\notag
\end{align}
where in the last identity we have used the fact that $\nu(q) = O(q^{-2+\ep})$.

We then, by inclusion-exclusion as in \cite[Section 4]{ShaTsi14}, have for an appropriate $Q$ (to be chosen later in terms of $X$),
\begin{align*}
N_{\operatorname{dom}}(U_e \cap \{a_{12}\ne 0 \}; X) &= \sideset{}{'}\sum_{q\geq 1} \mu(q) N_{\operatorname{dom}}(W_{q,e} \cap \{a_{12}\ne 0 \};X)\\
&=  \sideset{}{'}\sum_{1\leq q \leq Q} \mu(q) N_{\operatorname{dom}}(W_{q,e} \cap \{a_{12}\ne 0 \};X) \\
& \; \; \; +\;  \sideset{}{'}\sum_{q>Q} \mu(q) N_{\operatorname{dom}}(W_{q,e} \cap \{a_{12}\ne 0 \};X)\\
&=  \sideset{}{'}\sum_{1\leq q \leq Q} \mu(q) N^*(W_{q,e} \cap \{a_{12}\ne 0 \};X) \\
& \; \; \; +\;  \sideset{}{'}\sum_{1\leq q \leq Q} \mu(q) \left(N_{\operatorname{dom}}(W_{q,e} \cap \{a_{12}\ne 0 \};X)-N^*(W_{q,e} \cap \{a_{12}\ne 0 \};X)\right)\\
& \; \; \; +\;  \sideset{}{'}\sum_{q>Q} \mu(q) N_{\operatorname{dom}}(W_{q,e} \cap \{a_{12}\ne 0 \};X),
\end{align*}
where the sums are over square-free $q$ relatively prime to $e$.

By \cite[Lemma 3]{ShaTsi14}, we have $N_{\operatorname{dom}}(W_{q,e};X)=O(q^{-2+\ep} X)$ and we use this for the sum for $q>Q$.
We will use Equation~\eqref{e:Nstar} for the first $1\leq q \leq Q$ sum.
For the second $1\leq q \leq Q$ sum, note that each lattice point corresponding to a non-domain of discriminant $D$ is counted with coefficient $-\sum_{1\leq q \leq Q, q|D} \mu(q)$, which is $O(D^\ep)=O(X^{\ep})$, and by \cite[Equation (8)]{ShaTsi14} there are at most $O(X^{199/200+\ep})$ lattice points corresponding to non-domains.  (This step, or something similar, should be added to the proof in \cite{ShaTsi14}.)

As a result, as long as $Q=O( X^{1/80} e^{-1/2})$,
\begin{align*}
N_{\operatorname{dom}}(U_e \cap \{a_{12}\ne 0 \}; X)
&=\sideset{}{'}\sum_{q \geq 1} \mu(q) \nu(q) \delta_P \prod_{p\in P} (1-\nu(p))c_0X+ O(E_1) +O(E_2) +O(E_3),
\end{align*}
where
\begin{align*}
E_1&=  \sum_{1\leq q \leq Q} O(q^{\ep}  e  X^{39/40})\\
E_2 &=O(X^{199/200+\ep}) \\
E_3 &= \sum_{q>Q}q^{-2+\ep} X \\
E_4 & = \sum_{q >Q} \mu(q) \nu(q) \delta_P \prod_{p\in P} (1-\nu(p))c_0 X.
\end{align*}
These terms trivially admit the estimates
\begin{align*}
E_1&= O(Q^{1+\ep}  e  X^{39/40} )\\
E_2&=O(X^{199/200+\ep}  )\\
E_3&=O(Q^{-1+\ep} X)\\
E_4& = O(Q^{-1+\ep} X),
\end{align*}
where in the last estimate we have used the fact that $\nu(q) = O(q^{-2+\ep})$. 

We take $Q=X^{1/80}e^{-1/2}$, and have
\begin{align*}
N_{\operatorname{dom}}(U_e \cap \{a_{12}\ne 0 \}; X)
&= \sideset{}{'}\sum_{q \geq 1} \mu(q) \nu(q) \delta_P \prod_{p\in P} (1-\nu(p))c_0X+O(e^{1/2}X^{79/80+\ep}+X^{199/200+\ep})\\
&= \prod_p (1 - \nu(p)) \delta_P c_0X+O(e^{1/2}X^{79/80+\ep}+X^{199/200+\ep}). 
\end{align*}
From \cite[Lemma 11]{Bha10a}, we have that $ N_{\operatorname{dom}}( \{a_{12}=0 \}; X)=O(X^{39/40})$, and so
 $$
 N_{\operatorname{dom}}(U_e \cap \{a_{12}=0 \}; X)=O(X^{39/40}).
 $$
 It follows that
 $$
N_5(X;P)= N_{\operatorname{dom}}(U_e ; X)
=\prod_p (1 - \nu(p)) \delta_P c_0X+O(e^{1/2}X^{79/80+\ep}+X^{199/200+\ep})
 $$
 We thus conclude Theorem~\ref{thm_field_5_plus} holds with, $c_5 = \prod_p (1 - \nu(p))c_0$ when we only choose unramified splitting types.   As at the end of Theorem~\ref{thm_field_4_plus}, we can apply the result we have just proven and inclusion-exclusion to prove Theorem~\ref{thm_field_5_plus} in general.

\section{Application of the sieve}\label{sec_sieve_app}
 \subsection{Summary of the asymptotic inputs to the sieve}
We now turn to the application of the sieve lemma to degree $d$ field extensions of $\Q$. 
Note that when applying the sieve, it is crucial to have error terms with explicit dependence on local conditions (such as we have derived in Theorems \ref{thm_field_4} and \ref{thm_field_5}): without such an explicit dependence, we would not have quantitative control of the right-hand side of the key sieve inequality in Lemma \ref{lemma_sieve}, since we would not have an explicit bound for $R_p$ in terms of $p$.

Let $\Ascr$ and $\Ascr_p$ (for each rational prime $p$) be the sets so that $\# \Ascr = N_d(X)$ and $\# \Ascr_p = N_d(X;p)$ (or $\tilde{N}_4(X)$, $\tilde{N}_4(X;p)$ in the case of $d=4$).  With these definitions, the quantity $E(\Ascr; z, \frac{1}{2}M(z))$ treated in the sieve  (Lemma \ref{lemma_sieve}), which we will now denote by $E_d(\Ascr; z, \frac{1}{2}M(z))$, is the number of degree $d$ extensions $K$ of $\Q$ with $0<D_K \leq X$ (up to isomorphism, and non-$D_4$ when $d=4$) such that there are at most $\frac{1}{2}M(z)$ primes $p \leq z$ that split completely in $K$. 

We recall the collection $\Bscr_d(X;Y,M)$ of bad fields, as defined in (\ref{Bset_dfn'}).
We will think of $Y=z = (X/2)^{\del_0}$ for $\del_0>0$ to be chosen precisely later, and define $M(z)$ as in (\ref{M_dfn}). 
 In particular, the set of bad fields satisfies
\[ \#\Bscr_d(X;(X/2)^{\del_0}, \frac{1}{2}M((X/2)^{\del_0} )) 
 	=  E_d(\Ascr; (X/2)^{\del_0}, \frac{1}{2}M((X/2)^{\del_0})).
\]

We will need to apply the sieve separately to fields of each degree, since in several cases the count for $N_d(X;p)$ takes a somewhat different form, but in an effort to unify the presentation, we re-state the asymptotics we will assume in more general form. We write the results of Lemma \ref{lemma_field_2}, Theorem \ref{thm_field_3}, Theorems  \ref{thm_field_4} and \ref{thm_field_5} as follows.

{\bf Quadratic fields:} for $\del_e$ as in (\ref{del_p_2}), there is some $\sig_2  > 0$ and $0< \tau_2 \leq 1/2$ such that 
\begin{eqnarray*}
N_2(X) & = & c_2X  + O(X^{\tau_2+\ep}) \\
N_2(X;e) &=& \del_e c_2X  + O(e^{\sig_2}X^{\tau_2+\ep}) .
\end{eqnarray*}

{\bf Cubic fields:} for $\del_e, \del'_e$ as in (\ref{del_p_3}), there is some $ \sig_3  > 0$ and $0< \tau_3 <5/6$ such that 
\begin{eqnarray*}
N_3(X) & = & c_3X + c_3'X^{5/6}  + O(X^{\tau_3+\ep}) \\
N_3(X;e) &=& \del_e c_3X + \del_e'c_3'X^{5/6} + O(e^{\sig_3}X^{\tau_3+\ep}) .
\end{eqnarray*}

{\bf Non-$D_4$-quartic fields:} for $\del_e$ as in (\ref{del_p_4}), there is some $\sig_4  > 0$ and $0 < \tau_4 <1$ such that
\begin{eqnarray*}
\tilde{N}_4(X) & = & c_4 X + O(X^{\tau_4+ \ep}) \\
\tilde{N}_4(X;e) &=& \del_e c_4X  + O(e^{\sig_4} X^{\tau_4+\ep}) .
 \end{eqnarray*}

{\bf Quintic fields:} for $\del_e$ as in (\ref{del_p_5}), there is some $\sig_5 > 0$ and $0 <  \tau_5< 1$ as well as some $0 < \ga < 1$ such that
\begin{eqnarray*}
N_5(X) & = & c_5 X + O(X^{\ga + \ep}) \\
 N_5(X;e) &=& \del_e c_5X  +O(e^{\sig_5}X^{\tau_5}) +  O( X^{\ga+\ep}).
 \end{eqnarray*}

The main result of the sieve in this context is the following:
\begin{prop}\label{prop_main_3}
With the notation as above, we have
\[E_d(\Ascr; (X/2)^{\del_0}, \frac{1}{2}M((X/2)^{\del_0})) \ll  X^{1 - \del_0+ \ep},\]
for any $\del_0 $ such that 
\beq\label{del_choice_d}
 \del_0 \leq \begin{cases}
	\frac{1  - \tau_d}{1 + 2\sig_d} & \text{if $d=2,4$} \\
	\min\{ \frac{1  - \tau_d}{1 + 2\sig_d} , 1/4\} & \text{if $d=3$},\\
		\min\{\frac{1-\tau_d}{1+2\sig_d}, 1 - \ga\} & \text{if $d=5$}.
		\end{cases}
		\eeq
Moreover, for any such $\del_0$ there exist positive real constants $c_0(d) <c_1(d)<1$ and $X_{d} = X_d(\del_0) \geq 1$ such that for all $X \geq X_d$, 
\beq\label{M_bounds}
 c_0(d) \frac{(X/2)^{\del_0}}{\log (X/2)^{\del_0}}  \leq M((X/2)^{\del_0}) \leq c_1(d) \frac{(X/2)^{\del_0}}{\log (X/2)^{\del_0}}.
 \eeq
\end{prop}
The requirement that $X \geq X_d$ simply is a quantification of the requirement that $X$ be sufficiently large, and will be incorporated later simply by enlarging certain implicit constants.

Proposition \ref{prop_main_3} immediately provides the upper bound we require for the bad set $\Bscr_d(X;Y,M)$ defined in (\ref{Bset_dfn'}), with an appropriate choice of the parameters $Y,M$.  As there are $\pi(Y) = Y(\log Y)^{-1} + O(Y (\log Y)^{-2})$ primes $p \leq Y$, we could of course only expect at most $Y(\log Y)^{-1}$ primes $p\leq Y$ to split completely in any given field. Proposition \ref{prop_main_3} shows that, up to a constant factor, this is a reasonable expectation, in that the mean $M((X/2)^{\del_0})$ is approximately $\mu \pi ((X/2)^{\del_0})$, for some $\mu \in [c_0(d), c_1(d)]$; moreover Proposition \ref{prop_main_3} provides an upper bound for the number of fields with $D_K \leq X$ in which at most $\frac{1}{2}M((X/2)^{\del_0})$ primes $p \leq (X/2)^{\del_0}$ split completely.

We will prove Proposition \ref{prop_main_3} case by case.

\subsection{Sieve for quadratic fields}
For notational convenience, in this section we write $\sig,\tau$ for $\sig_2,\tau_2$. 
We compute that for any prime $p$,
\[R_p = \# \Ascr_p - \del_p \# \Ascr = 
	N_2(X;p)  - \del_p N_2(X)
	 = O(p^\sig X^{\tau + \ep}).\]
Similarly, for distinct primes $p,q$ 
\[R_{pq} = \# \Ascr_{pq} - \del_p \del_q \# \Ascr = 
	N_2(X;pq)  - \del_p \del_q N_2(X)
	 = O(p^\sig q^\sig X^{\tau + \ep}).\]
Thus since $\#\Ascr \gg X$,
\[  \frac{1}{\# \Ascr} \sum_{p | P(z)} |R_p|  \ll z^{1+\sigma}X^{\tau-1+\ep} , \qquad \frac{1}{\# \Ascr} \sum_{p,q | P(z)} |R_{pq}| \ll z^{2+2\sigma}X^{\tau-1+\ep}.
\]
We compute
\[
U(z)  = \sum_{p | P(z)} \del_p = \frac{1}{2} \sum_{p | P(z)} \frac{1}{1+p^{-1}},
\]
from which we deduce that
	 \beq\label{U_size_2}
\frac{1}{3} z(\log z)^{-1} + O(z (\log z)^{-2}) \leq U(z) \leq  	 \frac{1}{2} z(\log z)^{-1} + O(z(\log z)^{-2}).
\eeq
Indeed, letting $\ep_p = (1+p^{-1})^{-1}$, the upper bound follows directly from the prime number theorem and the fact that $0<\ep_p<1$, while the lower bound only requires noticing
\[ U(z) \geq \frac{1}{2} \sum_{p | P(z)} \ep_2 = \frac{1}{3}  \sum_{p | P(z)} 1 =  \frac{1}{3} z(\log z)^{-1} + O(z (\log z)^{-2}).\]
We may compute the mean as in (\ref{mean}):
\[ M(z)= U(z) +  \frac{1}{\# \Ascr} \sum_{p | P(z)} R_p  = U(z) + O(z^{1+\sig} X^{\tau-1 + \ep}).\]
Recalling (\ref{U_size_2}) and that $z=(X/2)^{\del_0}$ for a parameter $\del_0$ to be chosen later, we see the last error term will be $< \frac{1}{2}U(z)$ for sufficiently large $X$ as long as
\beq\label{del_cond_2}
\del_0< \frac{1 - \tau}{\sig}.
\eeq
Assuming this, for sufficiently large $X$ we have 
\[c_0 z(\log z)^{-1} \leq \frac{1}{2} U(z) \leq M(z) \leq \frac{3}{2} U(z) \leq  c_1 z(\log z)^{-1}\]
for absolute constants $0<c_0<c_1 \leq 1$.
We  apply Lemma \ref{lemma_sieve} to see that
\begin{eqnarray*}
 E_2(\Ascr; z, \frac{1}{2}M(z)) & \ll & \frac{X^{1+\ep}}{z^2} \left({z}+ z^{2+2\sig}X^{\tau-1}
 	 +z(z^{1+\sig}X^{\tau-1}) +(z^{1+\sig}X^{\tau-1})^2 \right)\\
	 &\ll& X^\ep( Xz^{-1} + z^{2\sig}X^{\tau}),
	 \end{eqnarray*}
still assuming (\ref{del_cond_2}).
Balancing the terms in the last expression above would set 
\beq\label{del_choice}
\del_0 = (1-\tau)/(1+2\sig),
\eeq
 which certainly satisfies (\ref{del_cond_2}); as a consequence, for any $\del_0 \leq (1-\tau)/(1+2\sig)$, we obtain 
\[ E_2(\Ascr; (X/2)^{\del_0}, \frac{1}{2}M((X/2)^{\del_0}))  \ll X^{1 - \del_0 + \ep} ,\]
which proves Proposition \ref{prop_main_3} in the case of quadratic fields.

\subsection{Sieve for cubic fields}
For notational convenience, in this section we write $\sig,\tau$ for $\sig_3,\tau_3$. We compute that
\[ R_p = \# \Ascr_p - \del_p \# \Ascr 
	= c_3'(\del_p'-\del_p)X^{5/6} + O(p^{\sig}X^{\tau +\ep})
		= O(p^{-1/3}X^{5/6} + p^{\sig}X^{\tau+\ep}).
\]
For distinct primes $p,q$,
\[  R_{pq} 
	= c_3'(\delta'_p\delta'_q-\delta_p\delta_q)X^{5/6} +O(p^{\sig}q^{\sig}X^{\tau+\ep})
= O(p^{-1/3}X^{5/6} +q^{-1/3}X^{5/6} + p^{\sig}q^{\sig}X^{\tau+\ep}).
\]
 Since $\# \Ascr \gg X$, we may compute that
 \begin{eqnarray*}
 \frac{1}{\# \Ascr} \sum_{p | P(z)} |R_p| 
	&\ll & z^{2/3}X^{-1/6} + z^{1+\sig}X^{\tau-1+\ep},\\
 \frac{1}{\# \Ascr}  \sum_{p,q| P(z)} |R_{p,q}|  
	 &\ll &  z^{5/3}X^{-1/6} + z^{2+2\sig}X^{\tau-1+\ep}.
	\end{eqnarray*}
Next, we note that
\[
U(z)  = \sum_{p | P(z)} \del_p = \frac{1}{6} \sum_{p | P(z)} \frac{1}{1+p^{-1}+p^{-2}} =  \frac{1}{6} \sum_{p | P(z)} e_p, 
\]
say.
From this we can deduce (as in the case of quadratic fields) that
	 \beq\label{U_size}
\frac{2}{21} z(\log z)^{-1} + O(z (\log z)^{-2}) \leq U(z) \leq  	 \frac{1}{6} z(\log z)^{-1} + O(z(\log z)^{-2}).
\eeq
Finally, we compute the mean
\[ M(z) 
 	= U(z) +  \frac{1}{\# \Ascr} \sum_{p | P(z)} R_p 
 		= U(z) + O(z^{2/3}X^{-1/6} + z^{1+\sig}X^{\tau-1+\ep}).
\]
Recalling (\ref{U_size}) and that $z=(X/2)^{\del_0}$ for a parameter $\del_0$ to be chosen later, we see the last error term will be $< \frac{1}{2}U(z)$ for sufficiently large $X$ as long as the analogue of (\ref{del_cond_2}) holds, in which case
\[ c_0 z(\log z)^{-1} \leq \frac{1}{2} U(z) \leq M(z) \leq \frac{3}{2} U(z) \leq  c_1 z(\log z)^{-1}.\]
for absolute constants $0<c_0 < c_1 \leq 1$.

We now apply Lemma \ref{lemma_sieve}, which shows that 
\begin{multline*}
 E_3(\Ascr; z, \frac{1}{2}M(z))  \ll  \frac{X^{1+\ep}}{z^2} \left({z}+ \{z^{5/3}X^{-1/6}+z^{2+2\sig}X^{\tau-1}\}  
  \right. \\
 	\left.
	+z\{z^{2/3}X^{-1/6}	+z^{1+\sig}X^{\tau-1}\}  +  \{z^{2/3}X^{-1/6}+z^{1+\sig}X^{\tau-1}\}^2 \right).
\end{multline*}
As long as $\del_0 \leq 1/4$, we have $ z^{5/3}X^{-1/6} \ll z$; after further simplification and still assuming the analogue of (\ref{del_cond_2}), we see that  
\[ E_3(\Ascr; z, \frac{1}{2}M(z))  \ll X^\ep( Xz^{-1} + z^{2\sig}X^{\tau}).\]

This is optimized by choosing $\del_0$ as in (\ref{del_choice}) as before, which satisfies (\ref{del_cond_2}).
 In particular, for any $\del_0 \leq \min \{ 1/4, (1-\tau)/(1+2\sig)\}$, we obtain 
\[ E_3(\Ascr; (X/2)^{\del_0}, \frac{1}{2}M((X/2)^{\del_0}))  \ll X^{1 - \del_0 + \ep} ,\]
which proves Proposition \ref{prop_main_3} in the case of cubic fields.

\subsection{Sieve for non-$D_4$-quartic fields}\label{sec_sieve_4}
The case of non-$D_4$-quartic fields is very similar to that for real quadratic fields, thus we only mention the highlights, with $\sig,\tau$ denoting $\sig_4,\tau_4$. We have 
\begin{eqnarray*}
R_p &=& \# \Ascr_p - \del_p \# \Ascr =  O(p^{\sig}X^{\tau +\ep}),\\
  R_{pq} &=&  \# \Ascr_{pq} - \del_p \del_q \# \Ascr   = O(p^\sig q^\sig X^{\tau+\ep}),\\
U(z)  &=& \sum_{p | P(z)} \del_p = \frac{1}{24} \sum_{p | P(z)} \frac{1}{1+p^{-1}+2p^{-2}+p^{-3}} .
\end{eqnarray*}
We deduce that 
	 \beq\label{U_size_4}
\frac{1}{3 \cdot 17} z(\log z)^{-1} + O(z (\log z)^{-2}) \leq U(z) \leq  	 \frac{1}{24} z(\log z)^{-1} + O(z(\log z)^{-2}).
\eeq
Next we compute the mean
\[ M(z)= U(z) +  \frac{1}{\# \Ascr} \sum_{p | P(z)} R_p  = U(z) + O(z^{1+\sig} X^{\tau-1 + \ep}).\]
Recalling (\ref{U_size_4}) and that $z=(X/2)^{\del_0},$
we see that as long as the analogous condition to (\ref{del_cond_2}) holds and $X$ is sufficiently large,
\[ c_0 z(\log z)^{-1} \leq \frac{1}{2} U(z) \leq M(z) \leq \frac{3}{2} U(z) \leq  c_1 z(\log z)^{-1}\]
for absolute  constants $0<c_0<c_1 \leq 1$.

We  apply Lemma \ref{lemma_sieve} to see that under the assumption (\ref{del_cond_2})
\[ E_4(\Ascr; z, \frac{1}{2}M(z))  \ll  \frac{X^{1+\ep}}{z^2} \left({z}+ z^{2+2\sig}X^{\tau-1}
 	 +z(z^{1+\sig}X^{\tau-1}) + (z^{1+\sig}X^{\tau-1})^2\right) \ll X^\ep( Xz^{-1} + z^{2\sig}X^{\tau}),\]
so that \[ E_4(\Ascr; (X/2)^{\del_0}, \frac{1}{2}M((X/2)^{\del_0}))  \ll X^{1 - \del_0 + \ep} \]
for any $\del_0 \leq (1-\tau)/(1+2\sig)$.

\subsection{Sieve for quintic fields}\label{sec_sieve_5}
Finally, we apply  the sieve to quintic fields, denoting $\sig_5,\tau_5$ by $\sig, \tau$.
We compute that for any $p = O(X^\rho)$, 
\[ R_p = \# \Ascr_p - \del_p \# \Ascr =  O(X^\ep(p^{\sig} X^{\tau} + X^{\gamma})).\]
For distinct primes $p,q$,
\[  R_{pq} =  \# \Ascr_{pq} - \del_p \del_q \# \Ascr   = O(X^\ep(p^\sig q^\sig X^\tau +  X^{\ga})).
  \]
We compute
\[
U(z)  = \sum_{p | P(z)} \del_p = \frac{1}{120} \sum_{p | P(z)} \frac{1}{1+p^{-1}+2p^{-2}+2p^{-3}+p^{-4}} , 
\]
from which we deduce that
\[\frac{2}{ 15 \cdot 37} z(\log z)^{-1} + O(z (\log z)^{-2}) \leq U(z) \leq  	\frac{1}{120}z(\log z)^{-1} + O(z(\log z)^{-2}).
\]
The mean may be expressed as 
\[ M(z)= U(z) +  \frac{1}{\# \Ascr} \sum_{p | P(z)} R_p  = U(z) + O(X^\ep(z^{1+\sig}X^{\tau-1}+   z X^{\ga-1 + \ep})).\]
The last term will be $< \frac{1}{2}U(z)$ for sufficiently large $X$ as long as $ \ga <1$ and the analogous condition to (\ref{del_cond_2}) holds. 
Assuming this, we have 
\[ c_0 z(\log z)^{-1} \leq \frac{1}{2} U(z) \leq M(z) \leq \frac{3}{2} U(z) \leq c_1  z(\log z)^{-1}\]
for absolute constants $0<c_0<c_1 \leq 1.$

We  apply Lemma \ref{lemma_sieve} to see that under the assumptions $\tau, \ga < 1$ and (\ref{del_cond_2}),
\[ E_5(\Ascr; z, \frac{1}{2}M(z))  \ll  \frac{X^{1+\ep}}{z^2} \left({z}+ z^2 X^{-1}(z^{2\sig}X^{\tau }+X^{\ga})
 	 +zX^{-1}(z^{\sig} X^\tau + X^{\ga})  + z^2X^{-2}(z^\sig X^\tau + X^{\ga})^2\right).\]
After simplification, this shows
 \[ E_5(\Ascr; z, \frac{1}{2}M(z))  \ll  \frac{X^{1+\ep}}{z^2} ({z} + z^2 X^{\ga-1}  + z^{2+2\sig} X^{\tau -1})
	\ll  X^\ep( Xz^{-1} + X^\ga+ z^{2\sig}X^{\tau}).\]
Assuming $z = (X/2)^{\del_0}$, we may conclude that for any $\del_0 \leq \min\{ (1-\tau)/(1+2\sig), 1 - \ga\}$, we have
\[ E_5(\Ascr; (X/2)^{\del_0}, \frac{1}{2}M((X/2)^{\del_0}))  \ll X^{1 - \del_0 + \ep} .\]
This completes the proof of Proposition \ref{prop_main_3}.

 \section{Proof of the main theorem and corollaries}\label{sec_main_theorems}
\subsection{Proof of Theorem \ref{thm_B}}

We now derive Theorem \ref{thm_B} from Proposition \ref{prop_main_3}.
By definition, if $M_1 \leq M_2$ then 
$\Bscr_d(X; Y, M_1) \subseteq \Bscr_d(X;Y, M_2)$. 
If $X$ is sufficiently large so that (\ref{M_bounds}) holds, say $X \geq X_d(\del)$, we may apply (\ref{M_bounds}) to write
\[
 \#\Bscr_d(X;(X/2)^{\del}, \frac{1}{2}c_0(d)  \frac{(X/2)^{\del}}{\log (X/2)^{\del}}) 
\leq  \#\Bscr_d(X;(X/2)^{\del}, \frac{1}{2}M((X/2)^{\del} )) = E_d(\Ascr; (X/2)^{\del}, \frac{1}{2}M((X/2)^{\del})).
\]
 We then apply Proposition \ref{prop_main_3} and deduce that for $X \geq X_d(\del)$,
 \[
  \#\Bscr_d(X;(X/2)^{\del}, \frac{1}{2}c_0(d)  \frac{(X/2)^{\del}}{\log (X/2)^{\del}}) \ll X^{1-\del+\ep}
 \]
for every $\ep>0$, and for $\del$ constrained by (\ref{del_choice_d}). 
 When we make the constraints in (\ref{del_choice_d}) precise by applying the results of Lemma \ref{lemma_field_2}, Theorem \ref{thm_field_3}, Theorems \ref{thm_field_4} and \ref{thm_field_5}, we obtain the parameters defined in (\ref{del_choice_d_specific}). For any $\del$ satisfying (\ref{del_choice_d_specific}), we may remove the explicit assumption that $X \geq X_d(\del)$ by including an appropriate implicit constant, so that
 \beq\label{B_bound_final}
  \# \Bscr_d(X;(X/2)^{\del}, \frac{1}{2}c_0(d)  \frac{(X/2)^{\del}}{\log (X/2)^{\del}}) \ll_{d,\del,\ep} X^{1  - \del +\ep}
  \eeq
for every $X \geq 1$ and every $\ep>0$.

\subsection{Proof of Theorem \ref{thm_main}}
To derive Theorem \ref{thm_main} from Theorem \ref{thm_B},  we proceed via a standard dyadic argument, which we now make precise. Let $\ep>0$ be fixed and for this $\ep$, let the implied constant in Theorem \ref{lemma_M_primes} be denoted by $C_0 = C_0(d,\ep)$, so that (\ref{EllVen_ineq}) becomes 
\beq\label{EllVen_ineq'}
 |\Cl_K[\ell]| \leq C_0 D_K^{1/2 + \ep}M^{-1}.
 \eeq
 Fix any $\del< \frac{1}{2\ell(d-1)}$. Then if $K$ is a degree $d$ extension of $\Q$ with $D_K \in (X,2X]$ that is not in 
the bad set $\Bscr_d^0(X;X^{\del}, \frac{1}{2} c_0(d) X^{\del}/ \log X^\del)$, 
we see from (\ref{EllVen_ineq'}) that 
\begin{eqnarray*}
 |\Cl_K[\ell]| &\leq& C_0 ( \frac{1}{2} c_0(d))^{-1}D_K^{1/2 + \ep}X^{-\del} \log (X^\del) \nonumber\\
		&\leq& C_0'D_K^{1/2 -\del + \ep} \log (D_K^{\del}),
\end{eqnarray*}
where it suffices to take $C_0' = C_0 2^{1+\del} ( \frac{1}{2} c_0(d))^{-1}$.
Now we assume that $X$ is sufficiently large, say $X \geq C(d,\ell,\ep)$, so that for all $\del< \frac{1}{2\ell(d-1)}$, and for all $D_K \in (X,2X]$, we have $\log (D_K^\del) \leq D_K^\ep$. Under this assumption we have 
\beq\label{desired_bound}  |\Cl_K[\ell]|  \leq C_0'D_K^{1/2 -\del + 2\ep}
\eeq
for all these fields not in $\Bscr_d^0(X;X^{\del}, \frac{1}{2} c_0(d) X^{\del}/ \log X^\del)$.

Let $F_{d,\ell}^0(X;\del,\ep)$ denote the collection of fields $K/\Q$ of degree $d$ with $X < D_K \leq 2X$ that fail the bound (\ref{desired_bound}); we may conclude that for any  $\del< \frac{1}{2\ell(d-1)}$ and for all $X \geq C(d,\ell,\ep)$,
\beq\label{FB}
F_{d,\ell}^0(X;\del,\ep) \subseteq \Bscr_d^0(X,X^{\del}, \frac{1}{2} c_0(d) X^{\del}/ \log X^\del).
\eeq
Now let  $F_{d,\ell}(X;\del,\ep)$ denote the collection of fields $K/\Q$ of degree $d$ with $0<D_K \leq X$ that fail the bound (\ref{desired_bound}); then 
\[ F_{d,\ell}(X;\del,\ep) \subseteq \Union_{0 \leq j \leq \lceil \log_2 X \rceil} F_d^0(2^{j};\del,\ep)
.\]

Set $j_0$ to be the smallest $j$ such that $2^{j_0} \geq C(d,\ell,\ep)$. 
Then for $j \leq j_0$,  we apply the trivial bound, $\#F_{d,\ell}^0(2^j,\del,\ep) \ll 2^j$. (This bound is only ``trivial'' in the sense that we know by (\ref{num_fields}) how to count fields of degree $d$ with $0<D_K \leq X$, for $d \leq 5$.)
For $j > j_0$ we apply (\ref{FB}) to write
\[  \Union_{j_0 < j \leq \lceil \log_2 X \rceil} F_{d,\ell}^0(2^{j};\del,\ep)	\subseteq \Union_{j_0 < j \leq \lceil \log_2 X \rceil} \Bscr_d^0(2^j,2^{j\del}, \frac{1}{2} c_0(d) 2^{j\del}/ \log 2^{j\del}).
	\]
Trivially enlarging each of the last sets to the non-dyadic version $\Bscr_d(2^{j+1},2^{j\del}, \frac{1}{2} c_0(d) 2^{j\del}/ \log 2^{j\del})$ and applying the result of Theorem \ref{thm_B} to each such set, we obtain
\beq\label{Fdl}
 \# F_{d,\ell}(X;\del,\ep) \ll C(d,\ell,\ep) + \sum_{j_0 < j \leq \lceil \log_2 X \rceil} 2^{j (1-\del + \ep')} \ll_{c,d,\ell,\ep,\ep'} X^{1 - \del+\ep'}.
 \eeq
which now holds (with a sufficiently large implicit constant) for all $X \geq 1$, for all $\ep'>0$  arbitrarily small,
and for all  $\del< \min\{ \frac{1}{2 \ell (d-1)},\del_0(d)\}$ where $\del_0(d)$ is defined as in (\ref{del_choice_d_specific}) in Theorem \ref{thm_B}.
For sufficiently large $\ell$, the first constraint on $\del$ is a stronger constraint than  the second.

To be precise, we now break down into cases depending on $d$.
For $d=2,$ Theorem \ref{thm_main}  is implied in the case $\ell=2$ by Gauss genus theory, and in the case $\ell=3$ by the known asymptotic (\ref{DavHei71_3}). For integers $\ell \geq 4$, Theorem \ref{thm_main} follows from (\ref{Fdl}), since ${1/2\ell} =  \min\{ 1/6, \frac{1}{2\ell}\}$ for $\ell \geq 4$. (Of course, for primes $\ell \geq 5$ and imaginary quadratic fields, Theorem \ref{thm_main} is implied by the stronger result (\ref{HB_Pie_result}), or indeed by an earlier result of Soundararajan \cite{Sou00} that at most one imaginary quadratic field $K$ with  $D_K \in [X,2X]$ can have $|\Cl_K[\ell]| \gg D_K^{1/2 - 1/2\ell + \ep}$; see also Corollary 2.2 of \cite{PieHB14a}.)
For $d=3$,  Theorem \ref{thm_main} is implied for $\ell=2$ by the known asymptotic (\ref{Bhargava_2}), and for $\ell=3$ by the stronger known result (\ref{cubic_3_del}). The cases $\ell \geq 4$ are implied by (\ref{Fdl}), since ${1/4\ell} =  \min\{ 2/25, \frac{1}{4\ell}\}$ for $\ell \geq 4$.
For $d=4$,  Theorem \ref{thm_main} follows from  (\ref{Fdl}) since $1/6\ell =  \min \{ \frac{1}{48}, \frac{1}{6 \ell} \}$ for $\ell \geq 8$; the remaining cases of $\ell \leq 7$ follow from the choice $\del_0 = 1/48$. 
Finally, for $d=5$, Theorem \ref{thm_main} similarly follows from (\ref{Fdl}) since $1/8\ell = \min \{ \frac{1}{200}, \frac{1}{8\ell}\}$ for  $\ell \geq 25$; the remaining cases of $\ell \leq 24$ follow from  the choice $\del_0 = 1/200$.

Corollaries \ref{cor_main_d} and \ref{cor_main_small} now follow from Theorem \ref{thm_main}, or can be derived directly from Theorem \ref{thm_B}, as already demonstrated in Section \ref{sec_anatomy_badfields}.

\section{Appendix: counting quadratic fields}\label{sec_real_quad}

In this appendix we prove the following result, from which Lemma \ref{lemma_field_2} may be deduced immediately.
 \begin{prop}\label{prop_quad_fields}
 Let $P$ be a finite set of primes. For each prime $p\in P$ we choose a splitting type at $p$ and assign a corresponding density as follows:  
 \begin{align*}
 \del_p & :=   \frac{1}{2}(1+ p^{-1})^{-1} & \textrm{ for } p = \pfrak_1 \pfrak_2 \\
 \del_p & :=  	 \frac{1}{2}(1+ p^{-1})^{-1}	& \textrm{ for }				 p = \pfrak_1\\
 \del_p & :=  	(p+1)^{-1}					& \textrm{ for } \text{p ramified.}
 \end{align*}
Let $e = \prod_{p \in P} p$ and $\del_e = \prod_{p \in P} \del_p$. Let $N_2^\pm(X;P)$ denote the number of real (respectively imaginary) quadratic extensions of $\Q$ with fundamental discriminant $|D_K| \leq X$ such that for each $p \in P$, the prime $p$ splits in the quadratic field with splitting type chosen for $p$ above. Then
 \beq\label{field_2_real_asympt}
 N_2^\pm(X;P) = \del_e \left( \frac{1}{3} +  \frac{1}{6}\right) \frac{1}{\zeta(2)} X+ O(e\sqrt{X}).
 \eeq
 \end{prop}
 We remark that in (\ref{field_2_real_asympt}), the first term is contributed by fundamental discriminants $\con 1 \modd{4}$ and the second by fundamental discriminants $\con 0 \modd{4}$.
We prove the proposition explicitly for $N_2^+(X;P)$, and omit the analogous argument for $N_2^-(X;P)$.
Upon combining the counts for real and imaginary fields, this implies Lemma \ref{lemma_field_2} as a special case.

The proof is a simple elaboration on the classical method for counting square-free integers $\leq X$. Recall that, for a fundamental discriminant $D$, a prime $p$ is ramified in $\Q(\sqrt{D})$ precisely when $p|D$; otherwise a prime $p \ndiv D$ splits in $\mathbb{Q}(\sqrt{D})$ if the Kronecker symbol $\Leg{D}{p}$ evaluates as $+1,$ and is non-split if $\Leg{D}{p}=-1$ (see e.g. \cite[Theorem 10.3, Chapter 16]{Hua}). Thus for each unramified $p \in P$ we assign $\ep_p \in \{-1,+1\}$ according to the specified splitting type of $p$. Let $P_0$ be the set of ramified primes in $P$ and set $P' = P\setminus P_0$; define $e_0 = \prod_{p \in P_0}p$ and $e' = \prod_{p \in P'} p$. Then we may write
\[N^{+}_2(X;P) = \# \{ \text{fundamental discriminants}  \; 1 \leq n \leq X : e_0|n, \; \Leg{n}{p} = \ep_p, \forall p \in P'\} .\]
We will find a count for this by sieving for fundamental discriminants (that is, elements that are free of odd squares) in the following two sets:
\begin{align*}
 \Acal^{(1)} & = \{ 1 \leq n \leq X : n \con 1 \modd{4},  e_0 |n,  \Leg{n}{p} =\ep_p, \forall p \in P' \} \\
 \Acal^{(0)} &=  \{ 1 \leq n \leq X:  n \con 8,12 \modd{16},  e_0|n, \Leg{n}{p} = \ep_p, \forall p \in P' \}.
 \end{align*}
More generally, fix a power $g$ and define for any $b \modd{2^g}$ the set
\[ \Acal =  \{ 1 \leq n \leq X:  n \con b \modd{2^g},  e_0|n, \Leg{n}{p} = \ep_p, \forall p \in P' \}.\]
 For each odd prime $q$ let $\Acal_q = \{ n \in \Acal : q^2 |n \}$. Note that certainly $\Acal_q$ is empty as soon as $q > \sqrt{X}$; we let $M$ be the index of the greatest prime $q_M \leq \sqrt{X}$. We will denote by $\overline{\Acal}_q$ the complement $\Acal \setminus \Acal_q$.
We will deduce Proposition \ref{prop_quad_fields} from the following lemma:
\begin{lemma}\label{lemma_inclusion}
Let $\Acal$ be as above, with $P = P_0 \union P'$ a set of odd primes. Then
\[  \bigcap_{q \; \mathrm{odd}} \overline{\Acal}_q  =  \frac{X}{3 \cdot 2^{g-2} \zeta(2)}  \prod_{p \in P'}\del_p \prod_{p \in P_0} \del_p + O(e\sqrt{X}),\]
with $\del_p$ as defined in Proposition \ref{prop_quad_fields}.
\end{lemma}
If the set $P$ specified in Proposition \ref{prop_quad_fields}  is a set of odd primes, then the proposition follows immediately from this lemma, by applying it to $\Acal^{(1)}$ with $g=2$, $b=1$ and then partitioning $\Acal^{(0)}$ into two disjoint sets with $g=4$ and $b=8$ or $12$, respectively, and applying the lemma to each.

If $2$ belongs to the set $P$ specified in Proposition \ref{prop_quad_fields}, then we consider separately the case when $2$ is specified to be ramified or unramified. If $2 \in P_0$ then $\Acal^{(1)}$ is empty. We already have $2|n$ for every $n \in \Acal^{(0)}$, so we set $P_{00} = P_0 \setminus \{2\}$ and apply Lemma \ref{lemma_inclusion} to $\Acal^{(0)}$ with $P = P_{00} \union P'$ (as before, separating $\Acal^{(0)}$ into two disjoint sets  and applying the lemma to each). We obtain 
\[ \bigcap_{q \; \text{odd}} \overline{\Acal}_q =  2 \cdot \frac{X}{3 \cdot 4 \zeta(2)}  \prod_{p \in P'}\del_p \prod_{p \in P_{00}} \del_p
	+O(e\sqrt{X}) 
	= \del_2 \cdot \frac{X}{2 \zeta(2)}  \prod_{p \in P'}\del_p \prod_{p \in P_{00}} \del_p +O(e\sqrt{X}) ,
	\]
with $\del_2=1/3$,	as claimed.
	
If $2 \in P'$ then $\Acal^{(0)}$ is empty.  We recall that for $p=2$ and $n \con 1 \modd{4}$, the Kronecker symbol $\Leg{n}{2}=+1$ if $n \con 1 \modd{8}$ and $-1$ if $n \con 5 \modd{8}$. Thus if $2 \in P'$, we set $P'' = P' \setminus \{2\}$ and $\Acal^{(1)}$ becomes
\[ \Acal^{(1)} = \{ 1 \leq n \leq X : n \con b \modd{8},  e_0 |n,  \Leg{n}{p} =\ep_p, \forall p \in P'' \} ,\]
with $b=1$ if  the original specification was $\ep_2=+1$ and $b=5$ if $\ep_2 =-1$. Applying Lemma \ref{lemma_inclusion}, we see that 
\[ \bigcap_{q \; \text{odd}} \overline{\Acal}_q  =  \frac{X}{3 \cdot 2 \zeta(2)}  \prod_{p \in P''}\del_p \prod_{p \in P_{0}} \del_p +O(e\sqrt{X}) 
	= \del_2 \cdot \frac{X}{2 \zeta(2)}  \prod_{p \in P''}\del_p \prod_{p \in P_{0}} \del_p +O(e\sqrt{X}) ,
	\]
	with $\del_2=1/3$, again as claimed. This proves Proposition \ref{prop_quad_fields}.

We now prove Lemma \ref{lemma_inclusion}. By the inclusion-exclusion principle, 
\beq\label{inclusion}
 \bigcap_{q \; \text{odd}} \overline{\Acal}_q = \sum_{m=0}^M (-1)^m  \sum_{q_1 < \cdots< q_m} | \Acal_{q_1} \intersect \cdots \intersect \Acal_{q_m} |,
 \eeq
in which for the $m=0$ term we sum the full set $|\Acal|$. 
A priori, any fixed term in (\ref{inclusion}) can be written as 
\[
  | \Acal_{q_1} \intersect \cdots \intersect \Acal_{q_m} |
	= \# \{ n \leq X : q_1^2 \cdots q_m^2 |n,  n \con b \modd{2^g}, e_0|n, \Leg{n}{p}=\ep_p, \forall p \in P'\}.
\]
Denote the set on the right-hand side by $S$, and let $Q:=\{q_1,\ldots, q_m\}$.
We first observe that if any $p \in P'$ belongs to $Q$ then the set $S$ must be empty. Thus we may reduce to considering the case in which $P'$ and $Q$ are disjoint, in which case we will prove that 
\beq\label{empty}
\# S = \frac{1}{2^g} \frac{X}{q_1^2 \cdots q_m^2} \frac{\gcd(q_1\cdots q_m,e_0)}{e_0} \prod_{p \in P'} \frac{1}{2}\left( \frac{p-1}{p} \right) + O(e').
 \eeq
First note that if a prime $p$ in $P_0$ belongs to $Q$ as well, then the condition $q_1^2 \cdots q_m^2 |n$  already specifies that $p$ is ramified. Thus upon defining $e_{00} = \prod_{p \in P_0 \setminus Q} p$, we may deduce that
	\beq\label{S_divide}
S= \{ k \leq X(q_1^2 \cdots q_m^2e_{00})^{-1}: k \con b(q_1^2 \cdots q_m^2e_{00})^{-1} \modd{2^g}, \Leg{k}{p} = \ep_p', \forall p \in P'\},
\eeq
where for each $p \in P'$ we have defined $\ep_p' = \ep_p \Leg{e_{00}}{p}.$

We note that for any integer $K \geq 1$ and any residue class $b$ modulo $2^g,$ the quantity
 \[\# \{ k \leq K : k \con b \modd{2^g}, \Leg{k}{p} = \ep_p', \forall p \in P' \}\] 
 may be expressed as
\begin{align*}
	& \sum_{\bstack{a \modd{e'}}{(a,e')=1}} \left( \prod_{p \in P'} \frac{1}{2}\left( 1 + \ep_p' \Leg{a}{p}\right) \right)\sum_{\tstack{k \leq K}{k \con a \modd{e'}}{k \con b \modd{2^g}}} 1 \nonumber \\
	& =  \frac{1}{2^{|P'|}}  \sum_{\bstack{a \modd{e'}}{(a,e')=1}}\prod_{p \in P'}\left( 1 + \ep_p' \Leg{a}{p}\right)  \left( \frac{K}{2^ge'} + O(1) \right) \nonumber \\
		 & =  \left( \prod_{p \in P'}\frac{p-1}{p} \right) \frac{K}{2^{g+|P'|}} + 0 + O(e');\label{compute_0}
		\end{align*}
	all the intermediate terms vanish	by orthogonality of characters.
Applying this to $S$ in (\ref{S_divide}), we obtain
\[ \# S= \frac{X}{q_1^2 \cdots q_m^2e_{00}} \frac{1}{2^{g + |P'|}} \prod_{p \in P'}\left( \frac{p-1}{p} \right) + O(e'),\]
proving (\ref{empty}).

Applying (\ref{empty}) to the inclusion-exclusion in (\ref{inclusion}) shows that
\[ \bigcap_{q \; \text{odd}} \overline{\Acal}_q  = \sum_{\bstack{d \leq \sqrt{X}}{(d,2e')=1}} \mu(d) \left(   \frac{X}{2^g d^2} \frac{\gcd(d,e_0)}{e_0} \prod_{p \in P' } \frac{1}{2}\left( \frac{p-1}{p} \right) +  O(e')\right).\]
The error term contributes $O(e\sqrt{X})$, while the main term contributes 
\[    X \frac{1}{2^g} \left( \prod_{p \in P'} \frac{1}{2}\left( \frac{p-1}{p} \right)\right) \sum_{\bstack{d=1}{(d,2e')=1}}^\infty \frac{\mu(d)}{d^2} \frac{\gcd(d,e_0)}{e_0} + O(X\sum_{d>\sqrt{X}} \frac{1}{d^2}) .
	 \]
	 Here the error term is $O(\sqrt{X}),$ with an implied constant which may be taken to be independent of $P$. 
	 
We now simplify the main term. We note that since $P$ consists of odd primes and  $(e_0,e')=1$, 	 upon setting $d=\del f$ with $\del = \gcd(d,e_0)$, we have
	\begin{eqnarray*}
	\sum_{\bstack{d=1}{(d,2e')=1}}^\infty \frac{\mu(d)}{d^2} \frac{\gcd(d,e_0)}{e_0}
	 & = &\sum_{\del | e_0} \sum_{\bstack{f=1}{(\del f,2e')=(f,e_0)=1}}^\infty \frac{\mu(\del f)}{\del^2 f^2} \frac{\del}{e_0}\\
	  &= &\frac{1}{e_0} \left( \sum_{\del | e_0} \frac{\mu(\del)}{\del} \right) \left(\sum_{\bstack{f=1}{(f,2e'e_0)=1}}^\infty \frac{\mu( f)}{ f^2} \right)
	  .
	  \end{eqnarray*}
	 The sum over $\del| e_0$ is a multiplicative function with respect to $e_0$. For $p$ prime we have 
	 \[ \sum_{\del | p} \frac{\mu(\del)}{\del}  = 1 - \frac{1}{p}
	 \]
	 and thus for $e_0$ square-free we may compute by multiplicativity that
	 \[ \frac{1}{e_0}  \sum_{\del | e_0} \frac{\mu(\del)}{\del}  = \prod_{p \in P_0} \frac{p-1}{p^2}.\]
	 
	 We next recall  that for any $\Re(s) >1$ and any distinct primes $q_1, \ldots, q_r$,
	 \[ \left(\prod_{i=1}^r(1 - \frac{1}{q_i^s})\zeta(s) \right)^{-1}= \prod_{p \not\in \{ q_1,\ldots, q_r\}} \left( 1 - \frac{1}{p^s} \right) = \sum_{\bstack{d=1}{(d,\prod q_i)=1}}^\infty \frac{\mu(d)}{d^s}.\]
Thus 
\[ \sum_{\bstack{f=1}{(f,2e'e_0)=1}}^\infty \frac{\mu( f)}{ f^2}
	 =(1 - \frac{1}{2^2})^{-1} \prod_{p \in P} (1-\frac{1}{p^2})^{-1} \frac{1}{\zeta(2)}.
\]
Assembling this all together, we see that
\[\bigcap_{q \; \text{odd}} \overline{\Acal}_q  = \frac{X}{2^g} (1 - \frac{1}{2^2})^{-1}  \frac{1}{\zeta(2)} \prod_{p \in P'} \left( \frac{p-1}{2p} \frac{1}{1 - \frac{1}{p^2}} \right) \prod_{p \in P_0} \left( \frac{p-1}{p^2}  \frac{1}{1 - \frac{1}{p^2}} \right) + O(e\sqrt{X}).
\]
This reduces to 
\[\bigcap_{q  \; \text{odd}} \overline{\Acal}_q  = \frac{X}{3 \cdot 2^{g-2} \zeta(2)}  \prod_{p \in P'}\frac{1}{2} \frac{1}{(1+p^{-1})}  \prod_{p \in P_0} \frac{1}{1+p} + O(e\sqrt{X}),
\]
proving Lemma \ref{lemma_inclusion}, with $\del_p$ as in Proposition \ref{prop_quad_fields}.

\section*{Acknowledgements}
The authors thank Paul Pollack, Arul Shankar, Frank Thorne, and Jacob Tsimerman for very helpful conversations and the referee for helpful comments on exposition.
We thank Frank Thorne in particular for a suggestion that helped improved the exponent in Theorem~\ref{thm_field_5_plus}.
Ellenberg is partially supported by National Science Foundation Grant DMS-1402620 and a Guggenheim Fellowship.
Pierce has been partially supported by National Science Foundation Grant DMS-1402121 and CAREER grant DMS-1652173.
Wood is partially supported by an American Institute of Mathematics Five-Year Fellowship, a Packard Fellowship for Science and Engineering, a Sloan Research Fellowship, and National Science Foundation Grant DMS-1301690.
\label{endofproposal}

\bibliographystyle{alpha}
\bibliography{NoThBibliography}

\end{document}

%% file: format.tex
\newtheorem{letterthm}{Theorem}

\newtheorem{thm}{Theorem}[section]
\newtheorem*{thm*}{Theorem}
\newtheorem{prop}[thm]{Proposition}
\newtheorem{lemma}[thm]{Lemma}
\newtheorem{cor}{Corollary}[thm]

\newcommand{\ep}{\varepsilon}

\newcommand{\con}{\equiv}

\newcommand{\ndiv}{\nmid}
\newcommand{\modd}[1]{\; ( \text{mod} \; #1)}
\newcommand{\bstack}[2]{#1 \atop #2}
\newcommand{\tstack}[3]{#1\atop {#2 \atop #3}}

\newcommand{\intersect}{\cap}

\newcommand{\Union}{\bigcup}
\newcommand{\union}{\cup}

\newcommand{\al}{\alpha}

\newcommand{\ga}{\gamma}
\newcommand{\del}{\delta}

\newcommand{\sig}{\sigma}

\newcommand{\Acal}{\mathcal{A}}

\newcommand{\beq}{\begin{equation}}
\newcommand{\eeq}{\end{equation}}